\documentclass[11pt]{article}
\usepackage{amssymb,amsmath,latexsym, verbatim}

\usepackage{tikz}

\allowdisplaybreaks

\usepackage{amssymb,amsmath,latexsym}
\usepackage{epsfig}
\usepackage{mathrsfs}
\usepackage{color}
\usepackage{graphicx}
\usepackage{epic,eepic,wrapfig,color, ifthen}
\usepackage{url}
\usepackage{amsthm}
\usepackage{amsfonts}
\usepackage{pmboxdraw}
\usepackage{verbatim}
\usepackage{color, mathtools, textcomp}
\usepackage{kotex}
\usepackage{cite}
\usepackage[normalem]{ulem}
\usepackage{enumerate}
\usepackage[bottom]{footmisc}
\usepackage{tikz}

\usepackage{hyperref}

\mathtoolsset{showonlyrefs}

\def\1{{\bf 1}}
\def\nn{\nonumber}

\def\R {{\mathbb R}}  
 
\def\E {{\mathbb E}}

\def\la{{\langle}}
\def\ra{{\rangle}}

\numberwithin{equation}{section}
\def\qed{{\hfill $\Box$ \bigskip}}

\def\BB{{\mathcal B}}

\def\R{{\mathbb R}}

\def\S{{\bf S}}

\def\P{{\mathbb P}}
\def\N{{\mathbb N}}

\def\eps{\varepsilon}

\def\wt{\widetilde}
\def\pf{\noindent{\bf Proof. }}

\def\beq{\begin{equation}}
\def\eeq{\end{equation}}

\def\la{\langle}
\def\ra{\rangle}

\theoremstyle{plain}
\newtheorem{thm}{Theorem}[section]
\newtheorem{lem}[thm]{Lemma}

\newtheorem{cor}[thm]{Corollary}
\newtheorem{remark}[thm]{Remark}
\newtheorem{prop}[thm]{Proposition}

\newtheorem{example}[thm]{Example}

\theoremstyle{remark}

\def\ud{\mathrm{d}}

\begin{document}

\newcommand\blfootnote[1]{
    \begingroup
    \renewcommand\thefootnote{}\footnote{#1}
    \addtocounter{footnote}{-1}
    \endgroup
}

\title{Dichotomy in the small-time asymptotics of spectral heat content for L\'evy processes}
	\author{Jaehun Lee\thanks{This work was supported by the National Research Foundation of Korea(NRF) grant funded by the Korea government(MSIT) (No. RS-2019-NR040050).} and Hyunchul Park\thanks{Research  supported in part by Research and Creative Project Award from SUNY New Paltz.} }
	\maketitle

\begin{abstract}
We establish a dichotomy in the small-time asymptotic behavior of the spectral heat content (SHC) for symmetric, but not necessarily isotropic, L\'evy processes whose L\'evy density satisfies a weak lower scaling condition near zero. This dichotomy is governed by whether the process has unbounded or bounded variation. In the unbounded variation case, the leading asymptotic behavior of the SHC is determined by the expected supremum of the process projected in the normal direction near the boundary. In contrast, for processes with bounded variation, the SHC decays linearly in time. Our main result, Theorem  \ref{thm:main}, extends and unifies key results from \cite{GPS19}, \cite{KP24}, and \cite{PS22}, covering a broader class of non-isotropic L\'evy processes and offering a streamlined proof.
\end{abstract}


	
\section{Introduction}
The spectral heat content (SHC) represents the total heat that remains in a domain $D\subset \R^{d}$ at time $t>0$ with the Dirichlet boundary condition and a unit initial condition.
The heat equation involves the Laplacian, which is the infinitesimal generator of Brownian motions.
If the Brownian motions are replaced by other L\'evy processes, one obtains the spectral heat content for those L\'evy processes.
In the past decade, there were extensive interests on SHC for jump processes (see \cite{KP23, KP24, KP24-3, P20, PS19, PS22, PX23} and references therein).
The most well-known example was investigated in \cite{PS22}, where the authors studied the small-time asymptotic behavior of rotationally invariant (isotropic) stable processes (SSP) on $\R^{d}$, $d\geq 2$. 
In their main result, they proved that there are three different regimes for the small-time asymptotic behavior of SHC for SSP, where the three regimes are determined by the stability index of SSP.

The purpose of this paper is to prove that there is a dichotomy in the small-time asymptotic behavior of SHC for a large class of L\'evy processes. We consider L\'evy processes with L\'evy density $J(x)$ that is comparable to $\frac{1}{|x|^{d}\psi(|x|)}$, where $\psi(x)$ satisfies so-called weak lower scaling condition near zero.
Within this class, we demonstrate that the small-time asymptotics of SHC exhibit two distinct regimes, determined by whether the underlying process X has unbounded or bounded variation.  
The condition that $X$ has unbounded variation is that either it has a non-trivial diffusion part ($A\neq 0$), or there are enough small jumps ($\int_{0+}\frac{dr}{\psi(r)}=\infty$). 
In the unbounded variation case, the leading contribution to the heat loss arises from particles near the boundary exiting the domain in the direction of the outward normal---this being the most efficient escape route. The asymptotic behavior of the SHC in this regime is given by an integral over the boundary involving the expected minimum of the supremum of the one-dimensional projection of the process in the normal direction and a fixed constant (e.g., 1), with respect to the surface measure.
In contrast, when $X$ has bounded variation ($A = 0$ and $\int_{0+} \frac{dr}{\psi(r)} < \infty$), particles located throughout the domain can still exit by time $t$ with linear probability in $t$. As a result, the SHC decays at rate $t$, and all regions of the domain contribute to the leading term. In this case, the leading coefficient is given by the so-called perimeter of $D$ with respect to the process $X$.
To illustrate the usefulness of the main theorem, we consider the following example (see Theorem \ref{thm:main}, Remark \ref{remark:main}, and Section \ref{section:example} for proof).
\begin{example}
{\rm Let $d \geq 2$ and $D \subset \R^d$ be a $C^{1,1}$ open set. 
Let $\Psi(\xi)=|\xi|^{\alpha}\log^{\beta}(1+|\xi|^{\gamma})$ be a characteristic exponent of a L\'evy process $X$ in $\R^d$ with $\alpha,\alpha+2\beta,\gamma\in (0,2)$. If $\alpha\in [1,2)$,
\[
\lim_{t\to 0}\frac{|D|-Q_{D}^{X}(t)}{\E[\overline{Y}_{t}\wedge 1]}=|\partial D|,
\]
where  $Y$ is (any) one-dimensional projection of X onto a one-dimensional subspace, $\overline Y_t := \sup_{0<s \le t} (Y_s - Y_0)$ is its supremum process, and $|\partial D|$ is the perimeter of $D$.
If $\alpha\in (0,1)$,
\[
\lim_{t\to 0}\frac{|D|-Q_{D}^{X}(t)}{t}=\text{Per}_{X}(D),
\]
where $\text{Per}_{X}(D)=\int_{D}\int_{D^c}J(y-x)dydx$ is the perimeter of $D$ with respect to $X$. } 
\end{example}

Let us highlight the contribution of this paper into the existing literature. 
Firstly, as mentioned above, we show that the dichotomy exists in the asymptotic behavior of SHC for L\'evy processes.
In particular, Theorem \ref{thm:main} shows that the three different regimes in the main result of \cite{PS22} is in fact two regimes (see Remark \ref{remark:main}). 
Secondly, our main theorem covers SHC for non-isotropic processes, whereas most known results for SHC were done for isotropic processes. 
Thirdly, when the process $X$ has unbounded variation, the main theorem, Theorem \ref{thm:main}, can cover the critical case when the one-dimensional projection is not integrable.
To the best of authors' knowledge, the only known example where the small-time asymptotic behavior of SHC is known for this class of L\'evy process is the Cauchy process (1-stable process) in \cite{PS22}. 
However, the proof in  \cite{PS22} relied on a highly sophisticated argument involving comparisons between subordinate killed and killed subordinate stable processes and required an upper bound for the transition density of the supremum process of the Cauchy process.
These detailed information for other processes is not available until now and there seems to be little hope that the argument in \cite{PS22} can be applied to obtain the asymptotic behavior of the SHC for the critical case.
A generalization of \cite{PS22} into more general processes than SSP was done in \cite{KP24}, but processes there have a regularly varying exponent with index in $(1,2]$ which does not cover the critical case mentioned above. 
Finally, the proof of Theorem \ref{thm:main}  is more robust than those in\cite{KP24} or \cite{PS22}, and it applies to a broader class of processes, including those with jump densities of variable order  (see Remark \ref{remark:main} and Section \ref{section:example} for concrete examples).

We outline the main ideas behind the proof of our main result. The first part of Theorem \ref{thm:main} is when the process has unbounded variation, and it follows a similar path as \cite{KP24}, where the authors studied SHC for L\'evy processes with the regularly varying characteristic exponent with an index in $(1,2]$, corresponding to the case where projections are integrable.
However, since we aim to provide a robust proof that is applicable for all processes that have unbounded variation, the main techniques employed differ significantly from those in \cite{KP24}.
The crucial ingredients for proofs are the near-diagonal lower bound for the heat kernel (Proposition \ref{prop:NDL}), the lower and upper bounds for the exit probabilities from balls (Lemma \ref{l:survival}), and corresponding results for the projections onto subspaces (Corollary \ref{c:sur}).
While these estimates and bounds are similar to those in \cite{BKKL, CKL}, we only assume the \textit{lower weak scaling condition near zero} on the jump density, whereas two-sided weak scaling condition was imposed in \cite{BKKL, CKL}. Additionally, the underlying processes may include a diffusion component, requiring a tailored derivation of the results.
Using these heat kernel estimates and upper and lower bounds for the exit probabilities from balls, we prove three key lemmas, Lemmas \ref{lemma:half-space}, \ref{lemma:inner ball}, and \ref{lemma:outer ball}, where we prove the small-time asymptotics of events of exiting from half-spaces, inner balls, and complement of outer balls of the underlying $C^{1,1}$ open set.

 The event of exiting from inner balls or complement of outer balls could arise in two different scenarios.
The main contribution occurs when the projected process with respect to the normal vector exits the domain, which is the most efficient way of exiting, and this contributes the main term.
On the other hand, there is an alternative way of exiting, which happens when the projection onto the orthogonal complement exits the domain. Since $C^{1,1}$ open set satisfies the $R$-ball condition in \eqref{e:ball}, the event can be controlled by those associated with the inner ball and the complement of outer ball.
Along the subspace orthogonal to the normal direction, the projected process has to move much longer distance, at least the square root of the distance of the most efficient path, and from estimates for the exit probabilities from balls in Corollary \ref{c:sur}, this is at most in $O(t)$, which is negligible compared to the main contribution from the most efficient exit. 
Finally, these estimates from Lemmas \ref{lemma:half-space}, \ref{lemma:inner ball}, and \ref{lemma:outer ball} are combined to prove Theorem \ref{thm:main} using Lemma \ref{lemma:div}, which is an integration-based version of the pointwise estimate of \cite[Lemma 5]{Berg}, specifically tailored to suit our needs.

For the second part of Theorem \ref{thm:main}, we assume that $X$ has bounded variation, in which case a related result is already available in  \cite[Theorem 3.4]{GPS19}. We show that under this assumption, the process $X$ does not hit the boundary of $D$ upon exiting the domain. This property allows us to apply  \cite[Theorem 3.4]{GPS19} to establish the desired asymptotic behavior.

The organization of this paper is as follows. 
Section \ref{section:preliminaries} is a preliminary section explaining definitions and the statement of the main result.  
In Section \ref{section:proofs}, we provide the proof for the main result, Theorem \ref{thm:main}. 
In Section \ref{section:example}, we illustrate examples that our main result can be applied to determine the small-time asymptotic behavior of SHC.
Section \ref{section:appendix} contains technical results on near diagonal heat kernel estimates and exit probabilities from balls.
In Subsection \ref{subsection:HKE}, we prove the near diagonal heat kernel estimates in Proposition \ref{prop:NDL} and use this to bound the exit probabilities in Lemma \ref{l:survival}. In Subsection \ref{subsection:projection}, we prove similar results hold true uniformly for projection processes onto subspaces, meaning the bounds in Lemma \ref{l:projection} and Proposition \ref{c:sur} on projection processes $X^{H}$ only depend on the original L\'evy process $X$ and remain uniform across all subspaces $H$.

In this paper, we use $c_i$ to denote constants whose values are unimportant and may change from one appearance to another.
The notation $\P_{x}$ stands for the law of the underlying processes started at $x\in \R^d$ and $\E_{x}$ stands for expectation with respect to $\P_{x}$. 
For simplicity, we use $\P$ and $\E$ for $\P_{0}$ and $\E_{0}$, respectively.
Sometimes, we write $x=(\tilde{x},u)$, where $\tilde{x}\in \R^{d-1}$ and $u\in \R$, and use $\P_{(\tilde{x},u)}$ and $\E_{(\tilde{x},u)}$ for $\P_{x}$ and $\E_{x}$, respectively. For a subspace $H \subset \R^d$ and $x \in H$, we use $B_H(x,r)$ to denote an open ball in $H$ centered at $x$ and radius $r>0$. When $H = \R^d$, we use $B(x,r) := B_{\R^d}(x,r)$ for simplicity.
	
\section{Preliminaries}\label{section:preliminaries}
Let $X=\{X_{t}\}_{t\geq 0}$ be a L\'evy process in $\R^d$ with $d\ge 2$.
It follows from the L\'evy-Khintchine theorem (\cite[Theorem 8.1]{Sato})

\[
\E[e^{i\langle \xi, X_{t} \rangle}]=e^{-t\Psi(\xi)}, \quad \xi\in\R^{d},
\]
with
\[
\Psi(\xi)=-i \la\xi, \gamma\ra + \la\xi, A\xi\ra -\int_{\R^{d}}\left(e^{i \la\xi, x\ra}-1-i \la\xi, x\ra 1_{\{|x|\le 1\}} \right)\nu(\ud x),
\]
where $\gamma\in\R^{d}$, $A$ is a symmetric, nonnegative-definite matrix, and $\nu$ is a measure satisfying $\nu(\{0\})=0$ and $\int_{\R^d}(1\wedge |x|^2)\nu(\ud x)<\infty$.
The law of $X$ is uniquely characterized by the characteristic exponent $\Psi(\xi)$, and the triplet $(A,\gamma,\nu)$ is called the L\'evy triplet of $X$. 
Throughout this paper, we assume that $X=\{X_t\}_{t \ge 0}$ is a symmetric L\'evy process corresponding to the triplet $(A,0,J(x)\ud x)$, where $A$ is non-degenerate or $A=0$. The key assumptions of this paper are as follows.
\begin{enumerate}[(a)]	
	
\item[(A)] 	\textbf{Lower weak scaling condition near zero}

There exist a constant $R_0 \in (0,\infty]$, an increasing function $\psi : (0,R_0] \to (0,\infty)$ and constants $0 < C_1 \le C_2$ such that 
	\begin{equation}\label{e:J}
	 \frac{C_1}{|x|^d \psi(|x|)} \le J(x) \le \frac{C_2}{|x|^d \psi(|x|)}, \qquad |x| \le R_0.
	 \end{equation} 
Moreover, the following lower weak scaling condition near zero on $\psi$ holds: there exist constants $\alpha \in (0,2]$ and $C_\psi \in (0,1]$ such that for any $0<r \le R < R_0$,
	\begin{equation}\label{e:psi}
		\frac{\psi(R)}{\psi(r)} \ge C_\psi \left( \frac{R}{r} \right)^\alpha.
	\end{equation} 
	We impose one of the following assumptions.
\item[(B1)] \textbf{Processes with unbounded variation}

The process $X$ has unbounded variation. It follows from \cite[Theorem 21.9]{Sato} and the assumption \eqref{e:J} that the process has unbounded variation if and only if the following conditions hold.
\begin{equation}\label{e:cond1}
	\int_0^{R_0} \frac{\ud r}{\psi(r)} = \infty \quad \mbox{or} \quad A \mbox{ is non-degenerate,}
\end{equation}
\item[(B2)] \textbf{Processes with bounded variation}

The process $X$ has unbounded variation. Again, it follows from \cite[Theorem 21.9]{Sato} and the assumption \eqref{e:J} that this is equivalent to the following.
\begin{equation}\label{e:cond2}
	\int_0^{R_0} \frac{\ud r}{\psi(r)} < \infty \quad \mbox{and} \quad A=0.
\end{equation}

\end{enumerate}
	We note that once \eqref{e:psi} holds for some $\alpha\in (0,2]$, it holds for any $\beta<\alpha$ with the same constant $C_{\psi}$.  
Also, we remark that when $\alpha >1$, the first part of \eqref{e:cond1} is a consequence of \eqref{e:psi} since $\frac{C_\psi R_0^\alpha}{\psi(R_0)} r^{-\alpha} \le \frac{1}{\psi(r)}$ for $r \in (0,R_0)$. 
However, the maximal value
$$ 
\overline{\alpha} = \overline{\alpha}(\psi) := \sup\{ \alpha \in (0,2] : \exists \, C_\psi \in (0,1] \,\, \mbox{such that} \,\, \eqref{e:psi} \,\, \mbox{holds with some } C_\psi \mbox{ and } \psi \}. 
$$
can be less or equal to $1$ (see Example \ref{ex:less than 1} (v)), even though $\int_0^{R_0} \frac{\ud r}{\psi(r)} = \infty$.

An open set $D \subset \R^d$ is said to be $C^{1,1}$ if there exist $R_1>0$ and $\Lambda>0$ such that for every $z \in \partial D$, there exist a $C^{1,1}$ function $\phi=\phi_z : \R^{d-1} \to \R$ satisfying $\phi(0) = 0$, $\nabla \phi(0) = (0,\dots,0)$, $\Vert \phi \Vert_\infty \le \Lambda$, $|\nabla \phi(x_1) - \nabla \phi(x_2)| \le \Lambda|x_1 - x_2|$ for $x_1,x_2 \in \R^{d-1}$, and an orthonormal coordinate system $CS_z$ of $y=(y_1,\dots,y_{d-1},y_d):=(\wt y,y_d)$ with origin at $z$ such that $D \cap B(z,R_1) = \{y = (\wt y, y_d) \mbox{ in } CS_z : y_d > \phi(\wt y) \} \cap B(z,R_1)$. The pair $(R_1, \Lambda)$ is called the $C^{1,1}$ characteristics of $D$. It is well known that any bounded $C^{1,1}$ open set $D$ with characteristics $(R_1,\Lambda)$ satisfies the following uniform interior and exterior ball-condition: There exists $R = R(d,R_1,\Lambda)>0$ such that for any $z \in \partial D$, there exist open balls $B_1$ and $B_2$ of radius $R$ such that
\begin{equation}\label{e:ball}
	B_1 \subset D, \,\, B_2 \subset \R^d \setminus \overline D \mbox{ and } \partial B_1 \cap \partial B_2 = \{z\}.
\end{equation}
We say that $D$ satisfies the \textit{$R$-ball condition} if \eqref{e:ball} holds with $R>0$. 

Let $D$ be an open set in $\R^{d}$, and let $\tau^X_{D}=\inf\{t>0:X_{t}\notin D\}$ be the first exit time of the L\'evy process $X$ from $D$. 
The spectral heat content $Q^X_{D}(t)$ of $X$ at time $t>0$ is defined by 
\[
Q^X_{D}(t)=\int_{D}\P_{x}(\tau_{D}^X>t)\ud x.
\]
As the survival probability $\P_{x}(\tau_{D}^X>t)$ is the unique solution for the (non-local) heat equation with the Dirichlet \textit{exterior condition} and unit initial condition, the spectral heat content $Q^X_{D}(t)$ represents the total heat that remains in $D$ at time $t>0$ when the initial temperature of $D$ is one and the temperature of $D^{c}$ is kept identically zero. 

A projected process $X^{\nu}=\{X^{\nu}_t\}_{t\geq 0}$ onto a line generated by $\nu$ is a one-dimensional L\'evy process defined by $X^\nu_t := \la X_t, \nu\ra$, where $\nu \in \S^{d-1} := \{ y \in \R^d : |y| = 1\}$. 
For a one-dimensional stochastic process $Y=\{Y_{t}\}_{t\geq0}$, we define the supremum process $\overline{Y}=\{\overline{Y}_t\}_{t\geq 0}$ by $\overline Y_t := \sup_{0<s \le t} (Y_s - Y_0)$.

We are ready to state the main theorem. 			
		\begin{thm}\label{thm:main}
{\rm
Let  $D$ be a bounded $C^{1,1}$ open set in $\R^d$, $d\geq 2$. Suppose that $X$ is a symmetric L\'evy process in $\R^d$ satisfying Condition \textbf{(A)}.
\begin{enumerate}[(1)]
\item Assume that $X$ has unbounded variation. Then,
	\begin{equation*}
		\lim_{t \to 0} \frac{|D| - Q_D^X(t)}{ \int_{\partial D} \E[\overline {X^{\nu(y)}}_t  \land 1] S(\ud y) } = 1,
	\end{equation*} 
where $\nu(y)$ is the outer unit normal vector at $y \in \partial D$ and $S(\ud y)$ is the surface measure of $\partial D$.
\item Assume that $X$ has bounded variation. Then,
	\begin{equation*}
		\lim_{t \to 0} \frac{|D| - Q_D^X(t)}{t } = \text{Per}_{X}(D),
	\end{equation*} 
where $\text{Per}_{X}(D)=\int_{D}\int_{D^c}J(y-x)dydx$ is the perimeter of $D$ with respect to $X$. 
\end{enumerate}
Moreover, when $X$ has unbounded variation, $t=o(\E[\overline {X^{\nu(y)}}_t  \land 1] )$ for all $y\in \partial D$.
}
\end{thm}
	
\begin{remark}\label{remark:main}
\begin{enumerate}
\rm{
\item 
If $X$ is isotropic, any one-dimensional projection $X^\nu$ is identical in law for all $\nu \in \S^{d-1}$. 
Therefore, if we further assume that X has unbounded variation, the following result holds:
	\[
		\lim_{t \to 0} \frac{|D| - Q_D^X(t)}{\E[\overline Y_t \land 1]} = |\partial D|,
	\]
where $Y$ is (any) one-dimensional projection of $X$ onto a one-dimensional subspace.
	\item Let us compare Theorem \ref{thm:main} with \cite[Theorem 1.1]{PS22}. Let $Y^{(\beta)}=\{Y^{(\beta)}_{t}\}_{t\geq 0}$ be an isotropic $\beta$-stable process on $\R$. 
	For any $\beta\in (0,2)$, $\E[\overline{Y^{(\beta)}}_{t}]=t^{1/\beta}\E[\overline{Y^{(\beta)}}_{1}]$ by the scaling property and $\E[\overline{Y^{(\beta)}}_{1}]<\infty$ if $\beta\in (1,2)$. 
	On the other hand, it follows from the scaling property
	\[
	\E[\overline{Y^{(\beta)}}_{t} \wedge 1]=\int_{0}^{1}\P(\overline{Y^{(\beta)}}_{t}>u)\ud u=t^{1/\beta}\int_{0}^{t^{-1/\beta}}\P(\overline{Y^{(\beta)}}_{1}>v)\ud v\sim t^{1/\beta}\E[\overline{Y^{(\beta)}}_{1}],
	\]
	where the last relation follows from \cite[Proposition 4.1]{Val16} and the notation $f(t)\sim g(t)$ means $\displaystyle\lim_{t\to 0}\frac{f(t)}{g(t)}=1$. 
	
When $Y^{(1)}$ is the standard Cauchy process on $\R$ whose characteristic exponent is $\Psi(\xi)=|\xi|$, we observe that 
	\begin{equation}\label{e:cauchy}
	\E[\overline{Y^{(1)}}_{t} \wedge 1]=\int_{0}^{1}\P(\overline{Y^{(1)}}_{t}>u)\ud u=t\int_{0}^{1/t}\P(\overline{Y^{(1)}}_{1}>v)\ud v\sim \frac{t\ln(1/t)}{\pi},
	\end{equation}
	where the last follows from \cite[Proposition 4.3-(i)]{Val16}. 
	Hence, Theorem \ref{thm:main} recovers \cite[Theorem 1.1]{PS22} and it provides us with the probabilistic meaning of the expression $\frac{t\ln(1/t)}{\pi}$, which is asymptotically equal to $\E[\overline{Y^{(1)}}_{t} \wedge 1]$. The expression $\frac{t\ln(1/t)}{\pi}$ was computed for Cauchy process (1-stable process) in one dimension using the density
of the supremum process of the Cauchy process in \cite{Val16} and it was proved in \cite{PS22} that the same expressions hold true for the Cauchy processes in all dimensions $d\geq 2$. However,  in \cite{PS22}, the authors did not provide the reason \textit{why} the expression $\frac{t\ln(1/t)}{\pi}$ appear in high dimensions. This was partly due to the fact that the asymptotic limit for isotropic $\beta$-stable processes was given by $\E[\overline{Y^{(\beta)}}_{t}]<\infty$, $\beta\in (1,2)$, whereas $\E[\overline{Y^{(1)}}_{t} ]=\infty$. 
In this paper, we show that the asymptotic limit of the spectral heat content for those processes is given by $\E[\overline{Y^{(\beta)}}_{t} \wedge 1]$, $\beta\in [1,2)$ including the case $\beta=1$.
This result clarifies that the two different regimes,  $\beta \in (1,2)$  and  $\beta = 1$ , in the main result of \cite{PS22} can, in fact, be unified into a single expression.

Finally, the second part of Theorem \ref{thm:main} covers (isotropic) $\alpha$-stable processes \cite[Theorem 1.1]{PS22} when $\alpha \in (0,1)$. For $\alpha$-stable processes, the condition \eqref{e:cond2} is equivalent to 
$\alpha \in (0,1)$. 
Therefore, Theorem \ref{thm:main} fully recovers \cite[Theorem 1.1]{PS22}.

	\item 
Theorem \ref{thm:main} recovers and generalizes \cite[Theorem 3.1]{KP24} which asserts that if the characteristic exponent $\Psi$ is in $\mathcal{R}_\beta(\infty)=\{f:\R_{+}\to \R_{+}: \displaystyle\lim_{x\to\infty}\frac{f(xy)}{f(x)}=y^{\beta} \text{ for all } y>0\}$ with $\beta\in (1,2)$,
\[ 
\lim_{t \to 0} \frac{|D| - Q^X_D(t)}{\Psi^{-1}(t)\E[\overline{Y^{(\beta)}}_1]} = |\partial D|, 
\] 
See Example \ref{e:4.1} for details. 
	\item 
	We remark that the denominator $\int_{\partial D}\E[\overline {X^{\nu(y)}}_t  \land 1]S(\ud y)$ in Theorem \ref{thm:main} cannot be replaced by $\int_{\partial D}\E[\overline {X^{\nu(y)}}_t ] S(\ud y)$. 
	The simplest case is the 1-stable process (Cauchy process). 
It is well-known that the Cauchy process is not integrable and $\E[\overline{Y^{(1)}}_t]=\infty$.
On the other hand, by \eqref{e:cauchy}, $\E[\overline{Y^{(1)} }_{t} \land 1] \sim \frac{t\ln(1/t)}{\pi} $.
	
	Another example is as follows. Let $\beta \in (1,2), \gamma \in (0,1)$ and $X$ be an isotropic and pure-jump L\'evy process corresponds to the L\'evy measure $J(\ud x) = |x|^{-d-\beta}\1_{\{|x| \le 1\}} + |x|^{-d-\gamma}\1_{\{|x| > 1\}}$. It follows from \cite[Theorem 25.3]{Sato} that $\E[\overline{X^\nu}_t ]=\infty$ for any $\nu\in\mathbf{S}^{d-1}$. }
	 \end{enumerate}
	\end{remark}

	\section{Proofs of Main Results}\label{section:proofs}
In this section, we provide the proof for Theorem \ref{thm:main}. Since the second part of Theorem \ref{thm:main} is relatively simple, we mainly focus on the first part of the theorem.
From now on, we always assume that $d \ge 2$, and the diffusion matrix $A$ of the symmetric L\'evy process $X$ in $\R^d$ satisfies that $A=0$ or $A$ is non-degenerate. Also, we assume Condition \textbf{(A)}. 
This is exactly the same as assumptions in Corollary \ref{c:sur}.
Since $X$ is a L\'evy process, $\int_{\R^{d}}(1\wedge |x|^2)J(x)\ud x<\infty$, and we set 
\[
C_3 := \int_{\R^d} \big( 1 \land \frac{|x|^2}{R_0^2}\big) J(x)\ud x < \infty.
\]
Note that we have
 \begin{equation}\label{e:C0-1}
J(B(0,R_0)^c) = \int_{B(0,R_0)^c} J(x)\ud x \le C_3
 	\end{equation}
 	and by \eqref{e:J},

 \begin{equation}\label{e:C0-2}
 \int_0^{R_0} \frac{r}{\psi(r)}\ud r = \omega_{d}^{-1}\int_{B(0,R_0)}\frac{|x|^{2-d}}{\psi(|x|)}\ud x   
 \leq \omega_{d}^{-1}R_{0}^{2}C_{1}^{-1}\int_{B(0,R_0)}\frac{|x|^2}{R_0^2} J(x) \ud x \leq \omega_{d}^{-1}R_{0}^{2}C_{1}^{-1}C_3,
\end{equation} 
where $\omega_{d}=\frac{2\pi^{d/2}}{\Gamma(d/2)}$ is the surface area of the unit ball in $\R^{d}$.

Following the definition in \cite[Equation (1.10)]{BKKL} (or \cite[Equation (2.20)]{CKL}), we define the \textit{scale function} $\phi:(0,R_0] \to (0,\infty)$ of $X$ by
\begin{equation}\label{e:phi}
	\phi(r):=  \frac{r^2}{\Vert A \Vert + 2\int_0^r \frac{s\ud s}{\psi(s)}}, \qquad r \le R_0.
\end{equation}

The following lemma is a tail estimate of the L\'evy measure in terms of $\psi$.
\begin{lem}\label{l:TJ}
	There exist constants $c_1,c_2>0$ such that for any $r \le R_0/2$,
	\begin{equation}\label{e:TJ}
		\frac{c_1}{\psi(2r)} \le J(B(0,r)^c) \le \frac{c_2}{\psi(r)}.
	\end{equation}
	The constants $c_1$ and $c_2$ depend on $\alpha$, $C_1$, $C_2$, $C_3$, $C_\psi$,  $d$ and $\psi(R_0)$.
\end{lem}
\pf 
We observe that by \eqref{e:J} and \eqref{e:psi}, for any $r \le R_0/2$,
\begin{align*}
&J((B(0,r)^c))=J\left(B(0,r)^c\cap B(0,R_0)\right) +J(B(0,R_0)^c)  \\
=&\int_{B(0,r)^c\cap B(0,R_0)}J(x)\ud x +J(B(0,R_0)^c)\\
\leq&\int_{B(0,r)^c\cap B(0,R_0)}\frac{C_2}{|x|^d\psi(|x|)}\ud x +J(B(0,R_0)^c)\\
	\leq& \omega_{d}C_2  \int_r^{R_0} \frac{\ud s}{s\psi(s)} + J(B(0,R_0)^c)  \le \frac{\omega_dC_2  r^\alpha}{ C_\psi \psi(r)} \int_r^{R_0} \frac{\ud s}{s^{1+\alpha}} + C_3 \\ 
	\leq& \frac{1}{\psi(r)} \big( \frac{\omega_{d}C_2}{\alpha C_\psi} + C_3\psi(R_0) \big).
\end{align*}

On the other hand, since $\psi$ is increasing, by \eqref{e:J} we have
\begin{align*}
	J((B(0,r)^c)) 
	\ge \int_{B(0,r)^c}\frac{C_1}{|x|^d\psi(|x|)}\ud x\geq C_1\omega_{d} \int_r^{R_0} \frac{\ud s}{s\psi(s)} \ge C_1\omega_{d} \int_r^{2r} \frac{\ud s}{s\psi(s)} \ge \frac{C_1\omega_d}{2\psi(2r)}. 
\end{align*}
\qed

	Recall that for $x \in \R^d$ and $\nu \in \S^{d-1}$, we define  $H(x,\nu) := \{ y \in \R^{d}: \la\nu, y\ra < \la \nu,x \ra \}$. 
	When $x = 0$ and $\nu=-e_{d}=(0,\cdots, -1)$, we write $H$ for the usual upper half-space $H(0,-e_d)=\{x=(x_1,\cdots, x_d)\in \R^{d}: x_{d}>0\}$.
	Also, for $\nu \in \S^{d-1}$ we define $X_{t}^\nu = \la X_t, \nu\ra$ for the one-dimensional projection process and $\displaystyle\overline{X^\nu}_t := \sup_{0<s \le t} (X^\nu_s - X^\nu_0)$ is the supremum process of $X^\nu$. Since $\overline{X^\nu}_t \ge 0$, we have for $b>0$,
	\begin{equation*}
	\E[\overline{X^\nu}_{t} \land b] = \int_0^\infty \P(\overline{X^\nu}_{t} \land b \ge r)\ud r =  \int_0^b \P( \overline{X^\nu}_{t} \ge r)\ud r. 
	\end{equation*}
The next lemma is analogue to \cite[Lemma 3.4]{KP24}. 
	\begin{lem}\label{lemma:half-space}
Assume that $X$ has unbounded variation.
For any $x \in \R^d$, $\nu \in \S^{d-1}$ and constants $a,b \in (0,\infty)$,
		$$ \lim_{t \to 0} \frac{\int_0^a \P_{x - r \nu}(\tau_{H(x,\nu)}^X \le t)\ud r}{\E[\overline{X^\nu}_{t} \land b]} = 1. $$
	Moreover, the convergence is uniform for all $x \in \R^d$ and $\nu \in \S^{d-1}$.
	\end{lem}
	\pf Without loss of generality, we may assume that $x = 0$ and $\nu = -e_d$. Let $Y_{t} := X^{-e_d}_{t}$.
Note that $\P_{(\wt 0, r)}(\tau_H^X \le t) = \P(\overline Y_t \ge r)$ and this implies
\[
\int_0^a \P_{(\wt 0, r)}(\tau_H^X \le t) \ud r = \int_0^a \P(\overline Y_t \ge r)\ud r. 
\]

Since $Y_0=0$ almost surely, 
\[
\{\overline{Y}_{t}\geq r\}=\{\tau^{Y}_{(-\infty,r)}\leq t\}\subset \{\tau^{Y}_{(-r,r)}\leq t\} \quad \text{a.s.}
\]
and 
\[
\{\tau^{Y}_{(-r,r)}\leq t\}\subset \{\tau^{Y}_{(-r,\infty)}\leq t\} \cup \{\tau^{Y}_{(-\infty,r)}\leq t\} \quad \text{a.s.}
\]
Hence, 
\begin{equation}\label{eqn:H ub}
\P(\overline{Y}_{t}\geq r)\leq \P(\tau^{Y}_{(-r,r)}\leq t),
\end{equation}
and by the symmetry of $Y$
\beq\label{eqn:H lb}
\P(\tau^{Y}_{(-r,r)}\leq t)\leq \P(\tau^{Y}_{(-r,\infty)}\leq t) +\P(\tau^{Y}_{(-\infty,r)}\leq t) =2\P(\tau^{Y}_{(-r,\infty)}\leq t)=2\P(\overline{Y}_{t}\geq r).
\eeq
Hence, from \eqref{eqn:H ub}  and \eqref{e:sur1}, we have for any $a,b>0$,
\begin{equation}\label{e:3.4.1}
	\left|\int_a^b	\P( \overline Y_t \ge r) \ud r\right| \le \left|\int_a^b \P(\tau_{(-r,r)}^Y \le t) \ud r  \right| \le c_1 t \left|\int_a^b \frac{1}{\phi(r)}\ud r\right| \le  \frac{c_1|a-b| t}{\phi(a \land b)} = c_2 t,
\end{equation}
where the function $\phi$ is defined in \eqref{e:phi}.

Since Condition \textbf{(B1)} hold, either $\int_0^{R_0} \frac{dr}{\psi(r)} = \infty$ or $A \neq 0$ is non-degenerate. If $\int_0^{R_0} \frac{dr}{\psi(r)} = \infty$, from \eqref{eqn:H lb} and \eqref{e:sur2}, we have
\begin{align*}
	\E[\overline Y_t \land b]& = \int_0^b \P(\overline Y_t \ge r)\ud r \ge \frac12\int_0^{b \land r_0} \P( \tau^Y_{(-r,r)} \le t)\ud r  \\ &\ge c_3t \int_{\frac{\psi^{-1}(t)}{4\sqrt d}}^{b \land r_0} \frac{\ud r}{\psi(4 \sqrt d r)} = \frac{c_3 t}{4\sqrt d} \int_{\psi^{-1}(t)}^{4\sqrt d(b \land r_0)} \frac{\ud r}{\psi(r)}. 
\end{align*}
Since $\displaystyle\lim_{t \to 0}\psi^{-1}(t) = 0$ and $\int_0^{R_0} \frac{\ud r}{\psi(r)} = \infty$, we conclude that 
\begin{equation}\label{e:3.4.2}
\lim_{t \to 0} \frac{t}{\E[\overline Y_t \land b]} = 0.
\end{equation}

By a similar way, if $A \neq 0$ is non-degenerate, by \eqref{eqn:H lb}, \eqref{e:sur3} and change of variable with $r = \sqrt{t r_1}$ we have 
\begin{align*}
	\E[\overline Y_t \land b] \ge c_4 \int_{t^{1/2}}^{b \land r_0} \exp \big(-c_5 \frac{r^2}{t} \big) \ud r = c_4 t^{1/2} \int_1^{\big(\frac{b \land r_0}{t}\big)^{1/2}} r_1^{-1/2}\exp(-c_5 r_1)\ud r_1 \ge c_6 t^{1/2},
\end{align*}
which also implies \eqref{e:3.4.2}.

Hence, it follows from \eqref{e:3.4.1} and \eqref{e:3.4.2}
	$$ \lim_{t \to 0}  \frac{|\int_a^b \P(\overline Y_t \ge r)\ud r|}{\int_0^b \P(\overline Y_t \ge r)\ud r} \le \lim_{t \to 0} \frac{ c_2 t}{\E[\overline Y_t \land b]} = 0. $$
	
Therefore,
	\begin{align*}
			\lim_{t \to 0} \frac{\int_0^a \P_{(\tilde{0},r)}(\tau_H^X \le t)\ud r}{\E[\overline Y_t \land b]} &= 1 - \lim_{t \to 0} \frac{\int_a^b \P(\overline Y_t \ge u)\ud u}{\int_0^b \P(\overline Y_t \ge u)\ud u} = 1.\end{align*}
Since the constants $c_2$ and $c_3$ are independent of $\nu \in \S^{d-1}$ and $x \in \R^{d}$, the convergence is uniform on $(x, \nu)$ and we reach the conclusion.
	\qed
	
\noindent The next lemma is analogue to \cite[Lemma 3.5]{KP24}. 
Let $r_0 > 0$ be the constant in Corollary \ref{c:sur}.
	\begin{lem}\label{lemma:inner ball}
Assume that $X$ has unbounded variation.
Let $R \in (0,r_0]$, $a \in (0,R/2)$ and $B := B(x,R)$. Then, for any $y \in \partial B$ we have
		\[ 
		\lim_{t \to 0} \frac{\int_0^a \P_{y-r\nu} (\tau_B^X \le t) \ud r}{\E [\overline{X^\nu}_{t} \land 1]} = 1, 
		\]
		where $\nu \in \S^{d-1}$ is a unit vector satisfying $y-R\nu =x$. 
		Moreover, the convergence is uniform for all $y \in \partial B$.
	\end{lem}
	\pf Without loss of generality, we may assume that $x=(\wt 0, R)\in \R^{d}$, $y=0$, $\nu = -e_d$, and $B = B((\wt 0,R),R)$. 
	Denote $H:= H_{-e_d}$ and $Y_t = X^{e_d}_t$.
	Clearly $B\subset H$ and this implies that $\tau^X_B \le \tau^X_H$. Hence, it follows from Lemma \ref{lemma:half-space},
	\[
	\liminf_{t \to 0} \frac{\int_0^a \P_{(\wt 0,r)}(\tau_B^X \le t)\ud r}{\E [\overline Y_t\wedge 1]} \ge \liminf_{t \to 0} \frac{\int_0^a \P_{(\wt 0,r)}(\tau_H^X \le t)\ud r}{\E [\overline Y_t\wedge 1]} = 1.
	\]

We now establish the opposite inequality, where the $\liminf$ is replaced by $\limsup$.
The following picture would help understand the proof. 
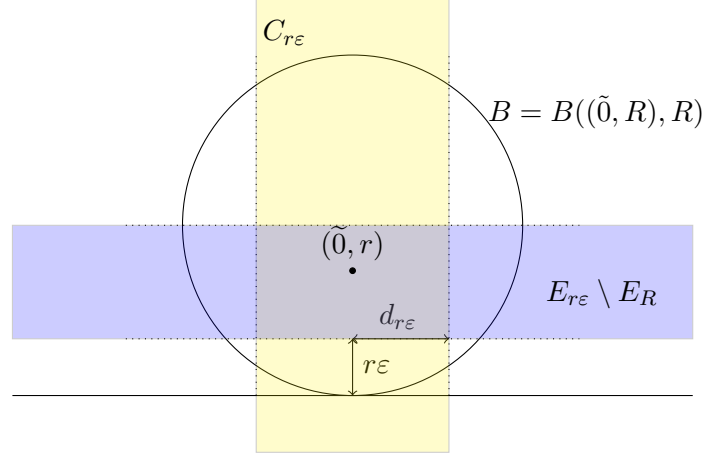
\begin{center}
\begin{tikzpicture}[scale=0.75]
\node at (7.3,5) {$B=B((\tilde{0},R),R)$};
\draw (-3,0)--(9,0);
\draw (3,3) circle [radius=3];
\draw [dotted] (-1,1)--(7,1);
\draw [dotted] (-1,3)--(7,3);
\draw [dotted] (1.3,0)--(1.3,6);
\draw [dotted] (4.7,0)--(4.7,6);
\draw [<->] (3,0)--(3,1);
\node [right] at (3,0.5) {$r \eps$};
\draw [<->] (3,1)--(4.7,1);
\node [above] at (3.8,1) {$d_{r\eps}$};
\draw [draw=black, fill=yellow, opacity=0.2]
       (1.3,-1) -- (1.3,7) -- (4.7,7) --(4.7,-1)-- cycle;
\node [above] at (1.8,6) {$C_{r\eps}$};
\draw [draw=black, fill=blue, opacity=0.2]
       (-3,1) -- (9,1) -- (9,3) --(-3,3)-- cycle;
\node at (7.4,1.8) {$E_{r\eps}\setminus E_{R}$};
\draw[fill] (3,2.2) circle [radius=0.05];
\node [above] at (3,2.2) {$(\wt 0, r)$};
\end{tikzpicture}
\end{center}
Fix $\eps \in (0,1)$. For $r \in (0,R/2)$ let $d_{r\eps} = (Rr\eps)^{1/2}$. Define $C_{r\eps} := \{ x= (\wt x , x_d) \in \R^d : |\wt x| \le d_{r\eps} \}$ and $E_{c} := \{ x= (\wt x, x_d) \in \R^d : c < x_d \}$. 
Since $(E_{r\eps} \setminus E_{R}) \cap C_{r\eps} \subset B$,
	\begin{equation}\label{eqn:decomposition}
	 \{ \tau_B^X \le t \} \subset \{ \tau_{C_{r\eps}}^X \le t \} \cup \{ \tau_{(E_{r\eps} \setminus E_R)}^X \le t\}. 
	 \end{equation}
	Note that for any subset $F,G\subset \R^{d}$ we have $\{\tau_{F\cap G}^{X}\leq t\}\subset \{\tau_{F}^{X}\leq t\} \cup \{\tau_{G}^{X}\leq t\}$, and this implies
\begin{equation}\label{eqn:EF}
	\P_z(\tau^X_{F \cap G} \le t) \le \P_{z}(\{\tau_{F}^{X}\leq t\} \cup \{\tau_{G}^{X}\leq t\}) \leq \P_z(\tau^X_F \le t) + \P_z(\tau^X_G \le t). 
\end{equation}
We write $X_{t}$ as $X_t=(X_t^{\tilde{H}}, Y_t)$, where $X_{t}^{\tilde{H}}=\pi_{\tilde{H}}(X_t)$ with $\tilde{H}=\{x=(\tilde{x}, 0)\in \R^d : \tilde{x}\in \R^{d-1}\}$. 
Hence, it follows from \eqref{eqn:decomposition} and \eqref{eqn:EF},
	\begin{align}\label{eqn:ball ub0}
\P_{(\wt 0,r)}(\tau_B^X \le t) &\le \P_{(\wt 0,r)}(\tau^X_{(E_{r\eps} \setminus E_R)} \le t) + \P_{(\wt 0,r)}(\tau^X_{C_{r\eps}} \le t) \nn \\
&\le \P_{(\wt 0,r)}(\tau_{E_{r\eps}}^X \le t) + \P_{(\wt 0,r)}(\tau^X_{E_R^{c}} \le t) +\P_{(\wt 0,r)}(\tau^X_{C_{r\eps}} \le t) \nn \\
&= \P_{-r}( \tau^Y_{(-\infty,-r\eps)} \le t ) + \P_{-r} (\tau^Y_{(-R,\infty)} \le t) + \P_{\wt 0}(\tau^{X^{\tilde{H}}}_{B(\wt 0,d_{r\eps})} \le t). 
	\end{align}
	
For the first term of \eqref{eqn:ball ub0}, 
	\begin{align}\label{eqn:ball ub1}
		\int_0^a \P_{-r}( \tau^Y_{(-\infty,-r \eps)} \le t ) \ud r 
		=\int_0^a \P( \overline Y_t \ge (1-\eps)r ) \ud r = \frac{1}{1-\eps} \int_0^{(1-\eps)a} \P(\overline Y_t \ge r)\ud r.
	\end{align}

For the second term of \eqref{eqn:ball ub0}, observe that $R-r\geq R/2$ since $r \in (0,R/2)$.
Hence, it follows from the symmetry of $Y_t$, \eqref{eqn:H ub} and \eqref{e:sur1},
	\[
	\P_{-r}(\tau^Y_{(-R,\infty)} \le t) =\P(\underline Y_t \le -(R - r))= \P(\overline Y_t \ge R - r) \le \P(\overline Y_t \ge R/2) \le \P(\tau^Y_{(-R/2,R/2)} \le t) \le \frac{c_1 t}{\phi(R/2)},\]
	where $\underline{Y}_{t}=\inf\{0\leq s\leq t : Y_{s} - Y_0\}$ is the infimum process of $Y$.
Thus,
	\beq\label{eqn:ball ub2} 
	\int_0^a \P_{-r}(\tau^Y_{(-R,\infty)} \le t) \ud r \le a \P(\overline Y_t \ge R/2) \le \frac{a c_1 t}{\phi(R/2)}.
	\eeq

Finally, we estimate the last term of \eqref{eqn:ball ub0}. Recall that  $a \in (0,R/2)$, and it follows from the change of variable $s = (Rr \eps)^{1/2}$
	\begin{align}\label{eqn:ball ub3}
\int_0^a \P_{\wt 0}(\tau^{X^{\tilde{H}}}_{B(\wt{0},d_{r\eps})} \le t)\ud r &= \int_0^a \P_{\wt 0}(\tau^{X^{\tilde{H}}}_{B(\wt 0, (Rr\eps)^{1/2})}  \le t  )\ud r 
= \frac{2}{R\eps} \int_0^{(Ra\eps)^{1/2}}  s\P_{\wt 0}(\tau^{X^{\tilde{H}}}_{B(\wt 0, s)} \le t )  \ud s.
	\end{align}
Since $s \le R \le r_0$, using \eqref{e:sur} for $X^{\tilde{H}}$ we have 
\begin{equation}\label{eqn:ball ub4}
	\P_{\wt 0}(\tau^{X^{\tilde{H}}}_{B(\wt 0,s)} \le t) \le c_2 \big( \frac{t}{ \psi(s/4)} + \exp(-\frac{c_3s}{t^{1/2}}) \big).
\end{equation}

Observe that by \eqref{e:C0-2} and by the change of variable $u=s/4$
$$ 
\int_0^{(Ra\eps)^{1/2}} \frac{s}{\psi(s/4)}\ud s \le \int_0^{r_0} \frac{s}{\psi(s/4)}\ud s = \int_0^{r_0/4} \frac{4u}{\psi(u)}\ud u  \le c_4,
$$
and by the change of variable $v = \frac{s}{t^{1/2}}$
$$ 
\int_0^\infty s \exp(-\frac{c_3 s}{t^{1/2}})\ud s = t\int_0^\infty v \exp(- c_3 v  )\ud v = c_5 t. 
$$
Combining above two inequalities with \eqref{eqn:ball ub3} and \eqref{eqn:ball ub4}, we have
\begin{equation}\label{eqn:ball ub5}
\int_0^a \P_{(\tilde{0},r)} (\tau^X_{C_{r\eps}} \le t)\ud r = \int_0^a \P_{\wt 0}(\tau^{X^{\tilde{H}}}_{B(\wt 0,d_{r\eps})} \le t) \ud r \le c_6 t,
	\end{equation}
where $c_6=c_{6}(\eps, R, \alpha, A, C_1, C_2, C_3, C_{\psi}, d)$.

Therefore, we conclude from \eqref{eqn:ball ub0}, \eqref{eqn:ball ub1}, \eqref{eqn:ball ub2}, and \eqref{eqn:ball ub5}
\begin{align} \label{eqn:ball ub6}
	\int_0^a \P_{(\wt 0,r)}(\tau_B^X \le t)\ud r  	 \le \frac{1}{1-\eps} \int_0^{a} \P_0(\overline Y_t \ge r)\ud r + c_6 t.
\end{align} 
Hence, it follows from \eqref{eqn:ball ub6}, Lemma \ref{lemma:half-space} and \eqref{e:3.4.2}
\[
\limsup_{t \to 0} \frac{\int_0^a \P_{(\wt 0,r)}(\tau_B^X \le t)\ud r}{\E [\overline Y_t\wedge 1]} \le \frac{1}{1-\eps}\limsup_{t \to 0} \frac{\int_0^a \P(\overline Y_t \ge r)\ud r}{\E[\overline Y_t \land 1]} + \limsup_{t \to 0} \frac{c_6 t}{\E[\overline Y_t \land 1]} \le \frac{1}{1-\eps}.
\]
Since $\eps>0$ is arbitrary, we conclude  that $\displaystyle\limsup_{t \to 0} \frac{\int_0^a \P_{(\wt 0,r)}(\tau_B^X \le t)\ud r}{\E [\overline Y_t\wedge 1]} \le 1$. 
Finally, we observe that the convergence is uniform on $y\in \partial B$ as the constant $c_6=c_{6}(\eps, R, \alpha, A, C_1, C_2,  C_3, C_{\psi}, d)$ does not depend on the directional vector $\nu$.
\qed

The next lemma is an analogue of \cite[Lemma 3.7]{KP24}. The proof is very similar to Lemma \ref{lemma:inner ball}, but we provide the details for the reader's convenience. 
\begin{lem}\label{lemma:outer ball}
Assume that $X$ has unbounded variation.
Let $R \in (0,r_0]$, $a \in (0,R/2)$, and $B := B(x,R)$. Then, for any $y \in \partial B$ we have
\[ 
\lim_{t \to 0} \frac{\int_0^a \P_{y-r\eta} (\tau^{X}_{B^c} \le t) \ud r}{\E [\overline{X^\eta}_{t} \land 1]} = 1,
\]
where $\eta\in \S^{d-1}$ is a unit vector satisfying $y+R\eta =x$.
Moreover, the convergence is uniform for all $y \in \partial B$.
\end{lem}
\pf Without loss of generality, we may assume that $y=0$, $B = B((\wt 0,-R),R)$, and $\eta = -e_d$. Let $Y_{t}= X^{-e_d}_{t}$ and $H = H_{-e_d}=\{x\in \R^{d}: x_{d}>0\}$. 
Clearly, $H\subset B^{c}$ and this implies  $\tau_H^X \le \tau^X_{B^c}$. Hence, it follows from Lemma \ref{lemma:half-space} we have
\[ 
\limsup_{t \to 0} \frac{\int_0^a \P_{y-r\eta}(\tau^X_{B^c} \le t)\ud r}{\E[\overline X_t^\eta\wedge 1]} \le \limsup_{t \to 0}  \frac{\int_0^a \P_{y-r\eta}(\tau^X_H \le t)\ud r}{\E[\overline X_t^\eta\wedge 1]} = 1. 
\]

Now we focus on establishing the opposite inequality, where the $\limsup$ is replaced by $\liminf$.
Recall that for $\eps \in (0,1)$ and $r \in (0,R/2)$, we have defined in Lemma \ref{lemma:inner ball} that $d_{r\eps} = (R r\eps)^{1/2}$, $C_{r\eps} = \{ x \in \R^d : |\wt x| \le d_{r\eps}   \}$ and $E_c:= \{ x \in \R^d : c < x_d  \}$. Note that $(E_{-R} \setminus E_{-r\eps}) \cap C_{r\eps} \subset B$ for any $\eps \in (0,1)$ and $r \in (0,R/2)$.
The following picture would help understand the proof. 
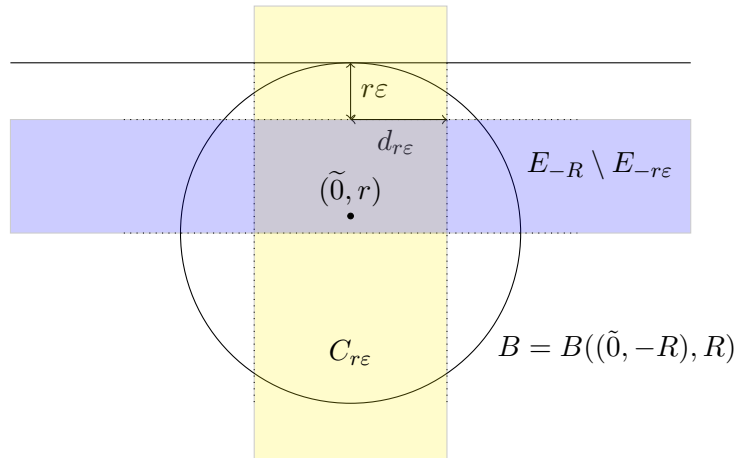
\begin{center}
\begin{tikzpicture}[scale=0.75]
\node at (7.7,-5) {$B=B((\tilde{0},-R),R)$};
\draw (-3,0)--(9,0);
\draw (3,-3) circle [radius=3];
\draw [dotted] (-1,-1)--(7,-1);
\draw [dotted] (-1,-3)--(7,-3);
\draw [dotted] (1.3,0)--(1.3,-6);
\draw [dotted] (4.7,0)--(4.7,-6);
\draw [<->] (3,0)--(3,-1);
\node [right] at (3,-0.5) {$r \eps$};
\draw [<->] (3,-1)--(4.7,-1);
\node [below] at (3.8,-1) {$d_{r\eps}$};
\draw [draw=black, fill=yellow, opacity=0.2]
       (1.3,1) -- (1.3,-7) -- (4.7,-7) --(4.7,1)-- cycle;
\node [above] at (3,-5.5) {$C_{r\eps}$};
\draw [draw=black, fill=blue, opacity=0.2]
       (-3,-1) -- (9,-1) -- (9,-3) --(-3,-3)-- cycle;
\node at (7.4,-1.8) {$E_{-R}\setminus E_{-r\eps}$};
\draw[fill] (3,-2.7) circle [radius=0.05];
\node [above] at (3,-2.7) {$(\wt 0, r)$};
\end{tikzpicture}
\end{center}

Note that
\begin{equation*}
 \{ \tau^X_{C_{r\eps}} > t \} \cap \{ \tau^X_{(E_{-R} \setminus E_{-r\eps})^c} \le t \} \subset \{ \tau^X_{B^c} \le t  \},
\end{equation*}
and this implies that 
$$
\{ \tau^X_{(E_{-R} \setminus E_{-r\eps})^c} \le t \} \subset \{ \tau^X_{B^c} \le t  \} \cup \{ \tau^X_{C_{r\eps}} \le t \}.
$$
Hence by \eqref{eqn:EF}, for any $\eps \in (0,R/2)$ and $y \in \R^d$ we have
\begin{equation*}
 \P_y(  \tau^X_{E_{-r\eps}} \le t) - \P_y(\tau^X_{E_{-R}} \le t) =  \P_y( \tau^X_{(E_{-R} \setminus E_{-r\eps })^c} \le t)  \le  \P_y(\tau^X_{B^c} \le t) + \P_y(\tau^X_{C_{r\eps}} \le t),
\end{equation*}
and this implies that for any $y\in \R^{d}$,
\begin{equation}\label{e:outer1}
\P_y(\tau^X_{B^c} \le t)\geq  \P_y(  \tau^X_{E_{-r \eps }} \le t) - \P_y(\tau^X_{E_{-R}} \le t) -\P_y(\tau^X_{C_{r\eps}} \le t).
\end{equation}

For the first term of the right hand side of \eqref{e:outer1},
\begin{align*}
	\int_0^a \P_{(\wt 0,r)}( \tau^X_{E_{-r \eps }} \le t ) \ud r =	
	\int_0^a \P( \overline Y_t \ge (1+\eps)r ) \ud r = \frac{1}{1+\eps} \int_0^{(1+\eps)a} \P( \overline Y_t \ge r)\ud r.
\end{align*}
For the second term of the right hand side of \eqref{e:outer1}, similarly to the argument in \eqref{eqn:ball ub2}, 
\begin{align*}
\int_{0}^{a}\P_{(\tilde{0},r)}(\tau^{X}_{E_{-R}}\leq t)\ud r
=\int_{0}^{a}\P(\overline{Y}_{t}\geq (r+R))\ud r \le a \P( \overline Y_t \ge R)
\leq c_{1}t,
\end{align*}
where the constant $c_1$ can be chosen independently of $\eta \in \S^{d-1}$.
Finally, it follows from \eqref{eqn:ball ub5} that the last term of the right hand side of \eqref{e:outer1} is bounded above by some constant times $t$, where the constant can be chosen independently of $\eta \in \S^{d-1}$.

Hence, it follows from \eqref{e:outer1}
\begin{align*}
 \int_0^a \P_{(\wt 0,r)}(\tau^X_{B^c} \le t) \ud r &\ge 	\int_0^a \P_{(\wt 0,r)}( \tau^X_{E_{-\eps r}} \le t ) \ud r - \int_0^a \P_{(\wt 0,r)}(\tau^X_{E_{-R}} \le t)\ud r - \int_0^a \P_{(\wt 0,r)}(\tau^X_{C_{r\eps}} \le t)\ud r \\ &\ge \frac{1}{1+\eps} \int_0^{(1+\eps)a} \P( \overline Y_t \ge r)\ud r - c_2t,
\end{align*}
where the constant $c_2$ can be chosen independently of $\eta \in \S^{d-1}$.
Using Lemma \ref{lemma:half-space} and \eqref{e:3.4.2}, we conclude that
$$ 
\liminf_{t \to 0} \frac{\int_0^a \P_{(\wt 0,r)}(\tau^X_{B^c} \le t)\ud r}{\E [\overline Y_t\wedge 1]} \ge \frac{1}{1+\eps} \liminf_{t \to 0} \frac{\int_0^{(1+\eps)a} \P( \overline Y_t \ge r) \ud r}{\E[\overline Y_t \wedge 1]} - \limsup_{t \to 0} \frac{c_2 t}{\E[\overline Y_t \wedge 1]} = \frac{1}{1+\eps}. 
$$
Since $\eps>0$ is arbitrary, this establishes the lower bound. 
\qed

We need the following lemma to prove Theorem \ref{thm:main}. The proof is similar in spirit to \cite[Lemma 5]{Berg}, but we need to modify it to fit our need.
For a subset $D\subset \R^{d}$ and $x\in D$, we define $\delta_{D}(x):=\inf\{|x-y|: y\in \partial D\}$ to be the Euclidean distance from $x$ to the boundary $\partial D$. For $a>0$, we define $D_{a}:=\{x\in D: \delta_{D}(x)>a\}$, which consists of points in $D$ whose distance to its boundary is at least $a>0$.
Since $D$ satisfies the $R$-ball condition \eqref{e:ball}, for any $x \in D \setminus D_a$, $a\leq R$, there exists unique $(y,r,\nu(x)) \in \partial D \times (0,a) \times \S^{d-1}$ such that $x = y -r\nu$, and $\nu(x)$ is the outer unit normal vector at $y \in \partial D$. 
Note that by the inner ball condition, for any $r_1,r_2 \in (0,a)$ and $y \in \partial D$, $\nu(y-r_1\tilde{\nu}) = \nu(y-r_2 \tilde{\nu})$, where $\tilde{\nu}\in \S^{d-1}$ is an outer unit normal vector in $\partial D$.
\begin{lem}\label{lemma:div}
Let $D$ be a bounded $C^{1,1}$ open set satisfying the $R$-ball condition.
Let $a \in (0,R/2)$. Then, for any bounded function $f:D \setminus D_a \to \R$,
\[
		\left( \frac{R-a}{R} \right)^{d-1} \le \frac{\int_{D \setminus D_a}  f(x)\ud x}{\int_0^a \int_{\partial D} f(y - s \nu(y)) S(\ud y)\ud s} \le \left(\frac{R}{R-a} \right)^{d-1},
\]
	where $S(\ud y)$ is the surface measure of $\partial D$, and $\nu(y)$ is the outer unit normal vector at $y\in \partial D$ such that $x=y-r\nu \in D\setminus D_{a}$.
\end{lem}
\pf 
Let $d(x) = dist(x,\partial D)$. It follows from the coarea formula that for any $f \in L^1(\overline{D \setminus D_a}, \ud x)$ 
\begin{equation*}
	\int_{D \setminus D_a} f(x)|\nabla d(x)|\ud x = \int_0^a \int_{\partial D_q} f(y)S_q(\ud y)\ud q,
\end{equation*}
where $S_q(\ud y)$ is surface measure of $\partial D_q$.
Since $\nabla d(x) = \nu(y)$, we have $|\nabla d(x)| = 1$ for all $x \in D \setminus D_a$. Thus, for $f \in L^1(\overline{D \setminus D_a}, \ud x)$,
\begin{equation}\label{e:coarea}
	\int_{D \setminus D_a} f(x)\ud x = \int_0^a \int_{\partial D_q} f(y)S_q(\ud y)\ud q.
\end{equation}
Since $D_s$ satisfies $(R-s)$ inner and outer ball conditions, by \cite[Equations (33) and (34)]{Berg} we have
\begin{equation}\label{e:div0}
	|{\rm div}\, \nu(x)| \le \frac{d-1}{R-s}, \qquad x \in \partial D_s. 
\end{equation}

From now on, we follow the proof of \cite[Lemma 5]{Berg}.  
Fix $0<q \le a$, $u : \partial D_q \to \R$ be a bounded function defined on $\partial D_{q}$, and let $\{u_{n}\}_{n=1}^{\infty}$ be a sequence of $C^{1}$ functions on $\partial D_{q}$ converging to $u$.
Define a sequence of $C^1$-function $g_n : D \setminus D_q \to \R$ by $g_n(y-s\nu) := u_n(y-q\nu)$, $y \in \partial D$ and $s \in (0,q)$. 
For any $0 \le q_1 < q_2 \le q$, using \eqref{e:coarea} and the divergence theorem for $g_n$ on $D_{q_1} \setminus D_{q_2}$, we obtain
\begin{equation}\label{e:div1}	
\int_{D_{q_1} \setminus D_{q_2}} \la\nabla g_n(x), \nu(x)\ra \ud x = -\int_{D_{q_1} \setminus D_{q_2}} g_n(x) {\rm div}\, \nu(x)\ud x + \int_{\partial D_{q_2}} g_n(y) S_{q_2}(\ud y) - \int_{\partial D_{q_1}} g_n(y)S_{q_1}(\ud y).
\end{equation}
Since $\la\nabla g_n, \nu\ra \equiv 0$ by the definition of $g_n$, we obtain from \eqref{e:coarea}, \eqref{e:div0}, and  \eqref{e:div1}
\begin{align}\label{e:div11}
	&\Bigg|\int_{\partial D_{q_1}} g_n(y)S_{q_1}(\ud y) - \int_{\partial D_{q_2}} g_n(y)S_{q_2}(\ud y)\Bigg| = \left| \int_{D_{q_1} \setminus D_{q_2}} g_n(x) {\rm div} \,\nu(x)\ud x \right| \nonumber \\ =& \left|\int_{q_1}^{q_2} \int_{\partial D_s} g_n(y){\rm div} \, (y)S_s(\ud y)\ud s \right| 
\le  \int_{q_1}^{q_2} \frac{d-1}{R-s} \Big( \int_{\partial D_s} g_n(y)S_s(\ud y)\Big)\ud s. 
\end{align}

For $s \in [0,q]$, let $h_n(s):= \int_{\partial D_{s}} g_n(y)S_s(\ud y)$. Then, by \eqref{e:div11},
\[ 
\left|h_n(q_1) - h_n(q_2)\right| \le  \left| \int_{q_1}^{q_2} \frac{d-1}{R-s} h_n(s)\ud s \right|. 
\]
Therefore, $|\frac{d}{\ud s} \log h_n(s)| \le \frac{d-1}{R-s}$, which implies
$$ 
\big(\frac{R-q}{R}\big)^{d-1}h_n(0) \le h_n(q) \le \big(\frac{R}{R-q}\big)^{d-1}h_n(0). 
$$
Since $h_n(0) = \int_{\partial D} g_n(y)S(dy) = \int_{\partial D} u_n(y+q\nu)S(dy)$ and $h_n(q) = \int_{\partial D_q} u_n(y)S_q(dy)$,
it follows from the bounded convergence theorem that for any $q \in [0,a]$,
\begin{equation}\label{e:div111}
\big( \frac{R-q}{R} \big)^{d-1} \int_{\partial D} u(y+q\nu)S(\ud y) \le \int_{\partial D_q} u(y)S_q(\ud y) \le \big( \frac{R}{R-q} \big)^{d-1}  \int_{\partial D} u(y+q\nu)S(\ud y). 
\end{equation}

Therefore, for any bounded function $f : D \setminus D_a \to \R$, by integrating both sides of \eqref{e:div111} and using \eqref{e:coarea} we have
$$ 
\int_0^a \big( \frac{R-a}{R} \big)^{d-1} \int_{\partial D} f(y+s\nu)S(\ud y)\ud s  \le  \int_{D \setminus D_a} f(x)\ud x = \int_0^a \int_{\partial D_q} f(y)S_q(\ud y)\ud s, 
$$
and
$$ 
\int_{D \setminus D_a} f(x)\ud x = \int_0^a \int_{\partial D_q} f(y)S_q(\ud y)ds \le \int_0^a \big( \frac{R}{R-a} \big)^{d-1} \int_{\partial D} f(y+s\nu)S(\ud y)\ud s. 
$$
This establishes the conclusion of this lemma. 
\qed

Now we are ready to prove Theorem \ref{thm:main}.

\noindent \textbf{Proof of Theorem \ref{thm:main}.} 
We focus on proving the first case when the underlying process $X$ has unbounded variation.
Since $D$ is $C^{1,1}$, by \eqref{e:ball} there exists $R>0$ such that $D$ satisfies the $R$-ball condition. Without loss of generality we assume that $R \le r_0$, where $r_0>0$ is the constant in Corollary \ref{c:sur}. 
For any $a < R/2$, we have $B(x,a) \subset D$ for all $x \in D_a$. Observe that for any $a < R/2$,
$$ |D| - Q_D^X(t) = \int_D \P_x(\tau_D^{X} \le t)\ud x = \int_{D_{a}} \P_x(\tau_D^X \le t)\ud x + \int_{D \setminus D_{a}} \P_x(\tau_D^X \le t)\ud x. $$
By \eqref{e:survival1}, for $x \in D_R$, we have $\P_x(\tau_D^{X} \le t) \le \P_x(\tau_{B(x,a)}^{X} \le t) \le \frac{c_1t}{\phi(a)}$. 
Therefore,
\begin{equation}\label{e:main1}
\int_{D_{a}} \P_x(\tau_D^X \le t)\ud x \le \frac{c_1t|D_a|}{\phi(a)} \le \frac{c_1t|D|}{\phi(a)}.
\end{equation}

Note that for any $\nu \in \S^{d-1}$, it follows from \eqref{e:3.4.2}  that $\displaystyle\lim_{t \to 0} \frac{t}{\E[\overline{X^\nu}_{t} \land 1]} = 0$.
Thus, it follows from the Fatou's lemma that for any $a < R/2$ we have
\begin{align}\label{e:main2}
&\limsup_{t \to 0} \frac{ \int_{D_a} \P_x(\tau_D^{X} \le t)\ud x}{  \int_{\partial D} \E[\overline{X^{\nu(y)}}_{t} \land 1] S(\ud y) } 
\leq \frac{c_{1}|D|/\phi(a)}{\displaystyle\liminf_{t\to 0} \int_{\partial D} \frac{ \E[\overline{X^{\nu(y)}}_{t} \land 1]}{t} S(\ud y)} \nonumber\\
\leq&  \frac{c_{1}|D|/\phi(a)}{ \int_{\partial D} \displaystyle \liminf_{t\to 0}\frac{\E[\overline{X^{\nu(y)}}_{t} \land 1]}{t} S(\ud y)}= 0,
\end{align}
where $\nu(y)$ is the outer unit normal vector at $y \in \partial D$.

 Next, we define $B_1(y) := B(y-R\nu(y), R) \subset D$ for $y \in \partial D$.
Note that by Lemma \ref{lemma:inner ball},
$$
\lim_{t \to 0} \frac{\int_0^a \P_{y-r\nu}(\tau_{B_1(y)}^X \le t)\ud r }{\E[\overline{X^{\nu}}_{t} \land 1] } = 1.
$$
Since the convergence is uniform on $y \in \partial D$ and $B_1(y) \subset D$ for all $y \in \partial D$ by the $R$-ball condition, 
$$ 
\limsup_{t \to 0} \frac{ \int_{\partial D} \int_0^a \P_{y-r \nu}(\tau_D^X \le t)\ud r S(\ud y)}{\int_{\partial D} \E[\overline{X^\nu}_{t} \land 1]S(\ud y)} \le \lim_{t \to 0} \frac{ \int_{\partial D} \int_0^a \P_{y-r \nu}(\tau_{B_1(y)}^X \le t)\ud rS(\ud y)}{\int_{\partial D} \E[\overline{X^\nu}_{t} \land 1]S(\ud y)} = 1. 
$$ 
Let $B_2(y) := B(y+R\nu(y), R)$. Then, $D \subset B_2(y)^c$ for all $y \in \partial D$ by the $R$-ball condition.
By repeating similar arguments as above with Lemma \ref{lemma:outer ball}, we obtain
$$ 
\liminf_{t \to 0} \frac{ \int_{\partial D} \int_0^a \P_{y-r \nu}(\tau_D^X \le t)\ud rS(\ud y)}{\int_{\partial D} \E[\overline{X^\nu}_{t} \land 1]S(\ud y)} \ge \lim_{t \to 0} \frac{ \int_{\partial D} \int_0^a \P_{y-r \nu}(\tau_{B_2(y)^c}^X \le t)\ud rS(\ud y)}{\int_{\partial D} \E[\overline{X^\nu}_{t} \land 1]S(\ud y)} = 1. 
$$
Therefore, by the above estimates we conclude that for any $a \in (0,R/2)$,
$$
\lim_{t \to 0}\frac{ \int_{\partial D} \int_0^a \P_{y-r \nu}(\tau_D^X \le t)\ud rS(\ud y)}{\int_{\partial D} \E[\overline{X^\nu}_{t} \land 1]S(\ud y)}  = 1. 
$$

Now observe that for $r \in (0,a)$ and $x = y-r\nu \in D \setminus D_a$, the function $f : D \setminus D_a \to [0,1]$ defined by $f(x) = \P_x(\tau_D^X \le t)$ is bounded function in $D \setminus D_a$. Applying Lemma \ref{lemma:div} for $f$ we obtain
$$ 
\left(\frac{R-a}{R} \right)^{d-1} \le  \frac{\int_{D \setminus D_a} \P_x(\tau_D^X \le t)\ud x }{\int_{\partial D} \int_0^a \P_{y-r \nu}(\tau_D^X \le t)\ud r S(\ud y) } \le \left(\frac{R}{R-a} \right)^{d-1}. 
$$
Combining this with \eqref{e:main2}, we obtain
$$ 
\left(\frac{R-a}{R} \right)^{d-1} \le \liminf_{t \to 0} \frac{ \int_{D} \P_{x}(\tau_D^X \le t)dx}{\int_{\partial D} \E[\overline{X^\nu}_{t} \land 1]S(\ud y)} \le \limsup_{t \to 0} \frac{ \int_D \P_{x}(\tau_D^X \le t)dx}{\int_{\partial D} \E[\overline{X^\nu}_{t} \land 1]S(\ud y)} \le \left(\frac{R}{R-a} \right)^{d-1}. 
$$
Since $a>0$ is arbitrary, we conclude from \eqref{e:main1} that 
\[
\lim_{t \to 0} \frac{|D| - Q_D^X(t)}{\int_{\partial D} \E[\overline{X^\nu}_{t} \land 1]S(\ud y)} = \lim_{t \to 0} \frac{\int_D \P_x(\tau_D^X \le t)\ud x}{\int_{\partial D} \E[\overline{X^\nu}_{t} \land 1]S(\ud y)} = 1. 
\]

Now we prove the second case when $X$ has bounded variation.
The main goal is to show that $X$ satisfies conditions for \cite[Theorem 1]{Sz00}. 
Note that $X$ has no Gaussian component, and the condition $\frac{C_1}{|x|^{d}\psi(|x|)}\leq J(x)\leq \frac{C_{2}}{|x|^{d}\psi(|x|)}$ for all $|x|\leq R_0$ ensures that the jumping kernel $J(x)$ is comparable for all directions. Hence, for any cone $V$ centered at the origin, there exists a constant $c(V)$ depending on the angle of the cone such that 
\[
\E_{0}\int_{0}^{\tau^{X}_{B(0,r)}}1_{V}(X_s)ds\geq c(V)\E_{0}\tau^{X}_{B(0,r)}.
\]
Since \cite[Theorem 1.3.1 and Lemma 1.3.2] {FOT} and the proof of \cite[Theorem 2.4]{CZ}, $X$ and $X^{D}$ have jointly continuous heat kernels. This implies that $X$ satisfies the strong Feller property and it is well-known that the semigroup strong Feller property implies the resolvent strong Feller property. 
It follows from Corollary \ref{c:5.3} that there exists a constant $c>0$ such that $\E_x \tau_{B(x,r)}^X \ge c\phi(r)$.
Hence,
\begin{align*}
&\int_{V\cap B(0,1)}\E_{0}\tau_{B(0,|y|)}J(dy)
\geq  c_1\int_{B(0,1)}\phi(|y|)J(dy)\\
\geq & c_2\int_{B(0,1)}\frac{\phi(|y|)}{\psi(|y|)|y|^d}dy
\geq c_{3}\int_{0}^{1}\frac{\phi(r)}{\psi(r)r}dr\\
=&\frac{c_{3}}{2}\int_{0}^{1}\frac{r}{\psi(r)}\left(\int_{0}^{r}\frac{s}{\psi(s)}ds\right)^{-1}dr\\
=&\frac{c_{3}}{2}\int_{0}^{1/2\phi(1)}\frac{dt}{t}=\infty,
\end{align*}
where we used the change of variable $t=\int_{0}^{r}\frac{s}{\psi(s)}ds$.
This shows that all conditions for \cite[Theorem 1]{Sz00} are met and $\P_{x}(X_{\tau_{D}}\in \partial D)=0$ for all $x\in D$. 
Now the conclusion follows immediately from \cite[Theorem 3.4]{GPS19}.
\qed

\section{Examples}\label{section:example}

In this section, we present examples of L\'evy processes for which our main result, Theorem \ref{thm:main}, yields the small-time asymptotic behavior of their spectral heat content.
\begin{example}\label{e:4.1}
	{\rm 
We verify that the first part of Theorem \ref{thm:main} generalizes the main result of \cite{KP24}.
Assume that $X$ is an isotropic L\'evy process with characteristic exponent $\Psi(x) = \Psi(|x|)$, and $\Psi \in \mathcal{R}_\beta(\infty)=\{f:\R_{+}\to \R_{+}: \displaystyle\lim_{x\to\infty}\frac{f(xy)}{f(x)}=y^{\beta} \text{ for all } y>0\}$ with $1<\beta < 2$.
Note that the function $\Psi : (0,\infty) \to (0,\infty)$ satisfy ${\rm WLSC}(\beta_1,\theta,c)$ and ${\rm WUSC}(\beta_2,\theta,C)$ for every $1<\beta_1< \beta < \beta_2<2$ (see \cite[Equations (17) and (18)]{BGR} for the definitions of WLSC and WUSC). Define a radial nondecreasing majorant of $\Psi$ by
$$ \Psi^*(r) := \sup_{|u| \le r} \Re \Psi(u) =  \sup_{s \le r} \Psi(s), $$
where the equality follows from the fact that $X$ is isotropic. (see \cite{G14} for details). In particular, since $\Psi$ satisfies ${\rm WLSC}(\beta_1,\theta,c)$, from the proof of \cite[Corollary 20]{BGR} the function $\Psi^*$ also satisfies ${\rm WLSC}(\beta_1,\theta',c')$ for some constants $\theta'>0$ and $c'>0$. 
It follows from  \cite[Equation (27)]{BGR}, \eqref{e:J} holds with the function $\psi(r) = \Psi^*(r^{-1})^{-1}$. Also, the condition ${\rm WLSC}(\beta_1,\theta',c')$ of $\Psi^*$ implies that the condition \eqref{e:psi} holds with index $\beta_1$.
Therefore, all conditions in Theorem \ref{thm:main} are satisfied.
}

{\rm Hence, the first or second part of Theorem \ref{thm:main} holds depending on whether $A\neq0$ or $\int_{0}^{R_0}\Psi(\frac{1}{r})dr=\infty$, or $\int_{0}^{R_0}\Psi(\frac{1}{r})dr<\infty$, respectively. 
This gives the following expression as a byproduct. If $A\neq0$ or $\int_{0}^{R_0}\Psi(\frac{1}{r})dr=\infty$, it follows from Theorem \ref{thm:main} and \cite[Theorem 3.1]{KP24} 
\[ 
\lim_{t \to 0} \frac{\E[\overline{X^\nu}_t \land 1]}{\Psi^{-1}(1/t)^{-1}\E[\overline{Y^{(\beta)}}_1] } = 1 \quad \text{for  } \beta\in (1,2),
\]
where $X^\nu$ is any one-dimensional projection of $X$, which are identical in law since $X$ is isotropic.		
}
\end{example}

\begin{example}\label{ex:less than 1}
{\rm
We illustrate some examples on the first part of Theorem \ref{thm:main}. These include processes with the critical characteristic exponent, truncated processes, and processes with variable order. 
}
\begin{itemize}
	\item[{\rm (i)}] {\rm(Critical characteristic exponent)} 
Assume that $X$ is a symmetric pure-jump L\'evy process whose L\'evy density satisfies \eqref{e:psi} with $\psi(r) = r \log(1+ \frac{1}{r})^{\beta}$ with $\beta>-1$. Then, $\int_{0}^{1}\frac{\ud r}{\psi(r)}=\infty$.
	In particular, if $\beta \in (-1,0)$, we have $\lim_{r \to 0} \frac{\psi(r)}{r} = 0$ and Theorem \ref{thm:main} entails L\'evy processes whose intensity of small jumps is higher than Cauchy process.    
	\item[{\rm (ii)}] {\rm (Truncated processes) Let $X$ be a L\'evy process with L\'evy measure $\nu^{X}(x)\ud x$ satisfying Conditions \textbf{(A)} and \textbf{(B1)}. Note that the conditions in Theorem \ref{thm:main} are unaffected when the jump density $\nu^{X}(x)\ud x$ is changed for $|x|>R_0$. 
Let $Y$ be a L\'evy process whose L\'evy density $\nu^{Y}(x)$ is $\nu^{Y}(x)=\nu^{X}(x)1_{\{|x|\leq R_0\}}$. Then, the first part of Theorem \ref{thm:main} still holds true for $Y$.
}
\item [{\rm (iii)}] {\rm (Jump diffusion)
Let $B$ be a Brownian motion in $\R^d$, and $Y$ be a symmetric pure jump process in $\R^d$ satisfying condition \textbf{(A)} which is independent of $B$. Let us consider a jump diffusion $X$ defined by $X_t = Y_t + B_t$. Since $X$ has unbounded variation, the process $X$ satisfies both Conditions \textbf{(A)} and \textbf{(B1)}. Therefore, we can apply the first part of Theorem \ref{thm:main} for $X$.}

	\item[{\rm (iv)}] {\rm (Processes with variable order) $\psi(r) = r^{\alpha(r)}$, $0<r<1$ where $1< \alpha_1 \le \alpha(r) \le \alpha_2 <2$. 
	Since $\psi(r)\geq \alpha_1>1$, the condition \eqref{e:cond1} is satisfied.
	Assume that $|\alpha'(r)| \le \frac{C}{r|\log r|}$. 
	Note that the function $f(r) = r \log r$ satisfies $\lim_{r \to 0} f(r) = 0$ and $f'(r) = 1+ \log r \le 0$ for $0<r \le e^{-1}$.
	Then, for any $0<r \le R \le \frac{1}{e}$,
	\begin{align*}
		&|\log R^{\alpha(R)-\alpha(r)}| = |(\alpha(R)-\alpha(r)) \log R| \\
		=& |\alpha'(r_0)(R-r) \log R| \leq |\alpha'(r_0)R \log R| \le \alpha'(r_0)r_0|\log r_0| \le c, 
	\end{align*}
	for some constant $c>0$.
	Hence, 
	\begin{align*}
		\frac{\psi(R)}{\psi(r)} = \frac{R^{\alpha(R)}}{r^{\alpha(r)}} \ge \big( \frac{R}{r} \big)^{\alpha(r)} \cdot R^{\alpha(R)-\alpha(r)} \ge e^{c} \big( \frac{R}{r} \big)^{\alpha_1}.
	\end{align*} 
	
	In particular, the function $\alpha(r) = \alpha_1 + r(\alpha_2 - \alpha_1)$ for $0<r<1$ satisfies the condition \eqref{e:psi} with $\alpha = \alpha_1$.
	}	
\end{itemize}
\end{example}

\begin{example}
{\rm 
	We illustrate an example of L\'evy process where the maximal of values of $\alpha$ that the condition \eqref{e:psi} holds is strictly less than 1, but $\int_{0}^{R_0}\frac{dr}{\psi(r)}=\infty$.
	Assume that $0<\alpha_1 < 1 < \alpha_2 < 2$. For $n \in \N$, define a strictly decreasing sequences $\{a_n\}_{n \ge 0}$ by $a_0 = 1$ and
	\[ 
	\int_{a_{n+1}}^{a_n} \frac{\ud r}{r| \log r|} = 1. 
	\]
	Note that $\inf_{n}\{a_{n}\}=0$. 
	Now we define $\alpha : (0,1) \to [\alpha_1, \alpha_2]$ as
	$$ \alpha(r) = \begin{cases} \alpha_2, &\quad r \in [a_{4n+1}, a_{4n}],\,\, n \in \N_0, \\
	\alpha_1 + (\alpha_2 - \alpha_1)\int_{a_{4n+2}}^r \frac{\ud r}{r|\log r|}, &\quad r \in [a_{4n+2},a_{4n+1}] \,\,n \in \N_0, \\
	\alpha_1, &\quad r \in [a_{4n+3}, a_{4n+2}],\,\, n \in \N_0, \\
			\alpha_2 + (\alpha_1 - \alpha_2)\int_{a_{4n+4}}^r \frac{\ud r}{r|\log r|}, &\quad r \in [a_{4n+4},a_{4n+3}] \,\, n \in \N_0. \\
	  \end{cases} $$
Then, we have $|\alpha(r)- \alpha(R)| \le \frac{R-r}{r|\log r|} \leq \frac{1}{r|\log r|}$ for $0<r<R<1$. Hence, the condition \eqref{e:psi} holds for $\psi(r) = r^{\alpha(r)}$with $\alpha= \alpha_1\in (0,1)$. 

\qquad Note that $\int_0^1 \frac{\ud r}{\psi(r)} \le \sum_{n=0}^\infty \int_{a_{4n+1}}^{a_{4n}} \frac{\ud r}{r^{\alpha_2}} \le \sum_{n=0}^\infty \int_{a_{4n+1}}^{a_{4n}} \frac{\ud r}{r |\log r|} = \infty$. 
Since $\psi(a_{4n+2}) = r^{\alpha_1}$, $\frac{\psi(a_2)}{\psi(a_{4n+2})} = \frac{a_2^{\alpha_1}}{a_{4n+2}^{\alpha_1}}\to \infty$ as $n\to \infty$. 
Hence, $\frac{\psi(R)}{\psi(r)}\geq c(\frac{R}{r})^{\alpha}$ for some fixed constant $c>0$ cannot hold for $0<r<R<1$ if $\alpha<\alpha_1$. 
This concludes that for $\psi(r)=r^{\alpha(r)}$ defined above, \eqref{e:cond1} holds and
\[ 
\sup\{ \alpha \in (0,2] : \exists \, C_\psi \in (0,1] \,\, \mbox{such that} \,\, \eqref{e:psi} \,\, \mbox{holds with } C_\psi \mbox{ and } \psi \} = \alpha_1 \in (0,1). 
\]}

\end{example}

\section{Appendix}\label{section:appendix}

	\subsection{Near diagonal estimates for Dirichlet heat kernel}\label{subsection:HKE}
In this section, we obtain the near diagonal lower bounds for the Dirichlet heat kernel and upper and lower bounds for survival probabilities for balls.
Recall that 
\begin{equation}
	\phi(r):=  \frac{r^2}{\Vert A \Vert + 2\int_0^r \frac{s\ud s}{\psi(s)}}, \qquad r \le R_0,
\end{equation}
is the scale function of $X$ defined in \eqref{e:phi}.
Note that
$\phi(r) \le \frac{r^2}{2\int_0^r \frac{s\ud s}{\psi(r)}} = \psi(r), \ r \le R_0. $
Note that by \cite[Theorem 1.3.1 and Lemma 1.3.2] {FOT} and the proof of \cite[Theorem 2.4]{CZ}, $X$ and $X^{D}$ have jointly continuous heat kernels (transition densities) $p_{X}(t,x,y)$ and $p_{X}^{D}(t,x,y)$, respectively. 
The following proposition is for the near-diagonal lower bound estimate for $p_X^{B(x_0,r)}(t,x,y)$.
\begin{prop}\label{prop:NDL}
	Let $X$ be a symmetric L\'evy process satisfying \eqref{e:J} and \eqref{e:psi}. Also, let us assume that $A=0$ or $A$ is non-degenerate. 
	Then, there exist $c>0$, $r_0>0$ and $\delta \in (0,1)$ such that for any $x_0 \in \R^d$,
	\begin{equation}\label{e:NDL}
		p_X^{B(x_0,r)}(t,x,y)  \ge c \phi^{-1}(t)^{-d}, \qquad r \le r_0, \,\, 0<t \le \phi(\delta r), \,\,x,y \in B(x_0,\delta\phi^{-1}(t)).
	\end{equation}
	The constants $c$, $r_0$ and $\delta$ depend only on $\alpha$, $A$, $C_1$, $C_2$, $C_3$, $C_\psi$, $d$ and $R_0$.
\end{prop}
\pf \textit{Step 1.} Consider extended function $\psi:(0,\infty) \to (0,\infty)$ as the function defined by
$$ \psi(r) := \begin{cases}    
	\psi(r), &\qquad r \le R_0, \\
	r^\alpha, &\qquad r > R_0,
\end{cases}  $$
and define $\phi(r)$ for $r>0$ as \eqref{e:phi}. Consider a symmetric pure-jump L\'evy process $X^1$ whose L\'evy measure $J_1(x)\ud x$ is given by
$$
J_1(x)=
\begin{cases}
	J(x), &\qquad |x| \le R_0, \\
	\frac{C_1}{|x|^{d+\alpha}} = \frac{C_1}{|x|^d \psi(|x|)}, &\qquad |x| > R_0.
\end{cases} 
$$
Then, $X^1$ satisfies \eqref{e:J} with the function $\psi$ and $R_0 = \infty$, and $\psi$ satisfies \eqref{e:psi} with the constants $\alpha$ and $C_\psi$. Also, from \cite[(2.3) and Lemma 2.4(2)]{BKKL} we can apply \cite[Theorem 3.11]{BKKL} for $X^1$. Thus, there exist constants $\delta_1 \in (0,1)$ and $c_1>0$ such that for any  $x_0 \in \R^d$,
\begin{equation}\label{e:ndl1}
	p_{X^1}^{B(x_0,r)}(t,x,y) \ge c_1 \phi^{-1}(t)^{-d}, \qquad r>0, \,\, 0<t \le \phi(\delta_1 r), \,\,x,y \in B(x_0,\delta_1 \phi^{-1}(t)).
\end{equation}
Moreover, the constants $c_1$ and $\delta$ depend only on $\alpha, C_1,C_2, C_3, C_\psi$ and $d$.
\\

\noindent \textit{Step 2.} Let $X^2_t:= A^{1/2}B_t + X^1_t:= \wt B_t + X^1_t$, where $(B_t)_{t \ge 0}$ is a Brownian motion in $\R^d$ which is independent of $X^1_t$. We claim that \eqref{e:NDL} holds for $X^2$. The case $A=0$ is trivial. Assume that $A \neq 0$ is symmetric and positive-definite matrix. 
By \cite[Theorem 1.1]{Zhang}, there is $c_2 = c_2(A, d)>0$ such that for any $x_0 \in \R^d$ and $C \ge 1$,
\begin{equation}\label{e:ndl2}
 p^{B(x_0,r)}_{\wt B}(t,x,y) \ge c_2 t^{-d/2}, \qquad r>0,\,\, 0<t \le Cr^2,\,\, x,y \in B(x_0,Ct^{1/2}). 
\end{equation}
Recall that by \eqref{e:C0-2} we have $2\int_0^r \frac{sds}{\psi(s)} \le 2\int_0^{R_0} \frac{sds}{\psi(s)} \le c_3$ for $r \in (0,R_0)$, where $c_3>0$ is a constant depends on $C_1$, $C_3$, $d$ and $R_0$. Thus, from  \eqref{e:C0-2} and \eqref{e:phi} we have 
\begin{equation}\label{e:square}
(c_3 + \Vert A \Vert)^{-1} r^2 \le \phi(r) \le \Vert A \Vert^{-1} r^2 \quad \mbox{for} \,\, 0<r \le R_0.
\end{equation}
Moreover, taking $r= \phi^{-1}(t)$ we have
\begin{equation}\label{e:inverse}
 \Vert A \Vert^{1/2} t^{1/2} \le \phi^{-1}(t) \le (c_3 + \Vert A \Vert)^{1/2} t^{1/2}, \quad t \in (0,\phi^{-1}(R_0)). 
 \end{equation}
Denote $\phi_1(r):= \frac{r^2}{\int_0^r \frac{s\ud s}{\psi(s)}} \ge \phi(r)$. Then, applying \eqref{e:ndl1} for $X^1$ we have for any $x_0 \in \R^d$,
\begin{equation}\label{e:ndl21}
 p_{X^1}^{B(x_0,r)}(t,x,y) \ge c_1\phi_1^{-1}(t)^{-d}, \qquad r > 0, \,\, 0<t \le \phi_1(\delta_1 r), \,\, x,y \in B(x_0,\delta_1 \phi_1^{-1}(t)). 
\end{equation}
Since $\phi(r) \le \phi_1(r)$ and $\phi_1^{-1}(t) \le \phi^{-1}(t)$, \eqref{e:ndl21} holds for $0<t \le \phi(\delta_1 r)$ and $x,y \in B(x_0,\delta_1 \phi_1^{-1}(t))$.
It follows from \eqref{e:ndl2} and \eqref{e:inverse}, there is $\delta_2 \in (0,1)$ such that
\begin{equation}\label{e:ndl22}
	p^{B(x_0,r)}_{\wt B}(t,x,y) \ge c_4 \phi^{-1}(t)^{-d}, \qquad 0<r \le R_0,\,\, 0<t \le \phi(\delta_2 r),\,\, x,y \in B(x_0, \delta_2 \phi_1^{-1}(t)). 
\end{equation}

Let $\delta_3 = \frac{\delta_1}{2} \land \frac{\delta_2}{2}$. Then, for any $x_0 \in \R^d$, $r \le R_0$, $0<t \le \phi(\delta_3 r)$ and $x,y \in B(x_0, \delta_3 \phi^{-1}(t))$, it follows from \eqref{e:ndl21} and \eqref{e:ndl22}, together with the fact that 
	$\{ \tau^{X^1}_{B(X^1_0,r/2)} > t \} \cap \{ \tau^{\wt B}_{B(\wt B_0,r/2)} > t \}  \subset  \{ \tau^{X^2}_{B(X^2_0,r)} > t \},$ 
\begin{align*}
&p^{B(x_0,r)}_{X^2}(t,x,y) \ge  \int_{\R^d} p^{B(x_0,r/2)}_{X^1}(t,x,z)p^{B(x_0,r/2)}_{\wt B}(t,z,y)\ud z \\
	&\ge \int_{B(x_0,\delta_2\phi_1^{-1}(t))} p^{B(x_0,r/2)}_{X^1}(t,x,z)p^{B(x_0,r/2)}_{\wt B}(t,z,y)\ud z \\
	&\ge \int_{B(x_0,\delta_2\phi_1^{-1}(t))} c_1\phi_1^{-1}(t)^{-d} \cdot c_4 \phi^{-1}(t)^{-d}  \ud z \ge   c_5 \phi^{-1}(t)^{-d}.\end{align*}
\\

\noindent \textit{Step 3.} Let $X^3$ be the $R_0$-truncated process of $X^2$, i.e., the L\'evy triplet of $X^3$ is given by $(A,0,J_1(\ud x)\1_{|x| \le R_0})$. 
In this time, we apply the L\'evy-Ito decomposition in \cite[Theorem 19.2]{Sato} to decompose $X^2$ into the truncated process $X^3$ and a compound Poisson process whose jump sizes are larger than $R_0$.
Define
\[ 
j(C,\omega) = \#\{s : (s,X^2_s(\omega) - X^2_{s-}(\omega)) \in C \cap \big((0,\infty) \times B(x,R_0)^c \big) \}, \qquad C \in \BB((0,\infty) \times (\R^d\setminus\{0\})). 
\]
Then, $j(C)$ is a Poisson random measure with the intensity measure $\nu$ with $\nu((0,t] \times E)):= tJ_1(E \cap B(x,R_0)^c)$ for $E \in \BB(\R^d\setminus\{0\})$. Hence,
\[ 
X^2_t = Y_t + \int_{\R^d} x j(\ud s \ud x):= Y_t + Z_t,
\]
where $Y$ and $Z$ are independent, and $Y$ and $X^3$ are identical in law.

For any $x_0 \in \R^d$, $r \le R_0/2$, $x \in B(x_0,r)$, and a Borel set $F \subset B(x_0,r)$, we have
$$    \P_{x}(\{ Y^{B(x_0,r)}_t \in F \} \cap \{Z_t = 0 \} )  = \P_{x}(X^{2,B(x_0,r)}_t \in F),  $$
where $X^{2,B(x_0,r)}$ is the killed process of $X^2$ upon $B(x_0,r)$. 
Indeed, since $r \le \frac{R_0}{2}$, we observe that $X_t^{B(x_0,r)} = \partial$ if there exists a jump whose size is larger than $R_0$. Therefore, since $X^3$ and $Y$ are identical in law, for any $x_0 \in \R^d$, $r \le r_0$,  and $0<t \le \phi(\delta_3 r)$ we have
$$ \P_{x}(X^{3,B(x_0,r)}_t \in F) = \P_{x}(Y^{B(x_0,r)}_t \in F) \ge \P_x(X^{2,B(x_0,r)}_t \in F), \quad x \in B(x_0,r), \,\, F \subset B(x_0,r). $$
Thus, for any $x_0 \in \R^d$, $r \le r_0$, and $0<t \le \phi(\delta_3 r)$, 
\begin{equation}\label{e:1234}
p_{X^3}^{B(x_0,r)}(t,x,y) \ge p_{X^2}^{B(x_0,r)}(t,x,y) \ge c_6 \phi^{-1}(t)^{-d}, \quad x,y \in B(x_0,\delta_3 \phi^{-1}(t)).
\end{equation}
\\

\noindent \textit{Step 4.} Since $X^3$ is also $R_0$-truncated process of $X$, by a similar way as before, we obtain
$$ \P_{x}(\{ X^{3,B(x_0,r)}_t \in F \} \cap \{\wt Z_t = 0 \} )  = \P_{x}(X^{B(x_0,r)}_t \in F), $$
where $\wt Z$ corresponds to a compound Poisson process independent of $X^3$ with the intensity measure $\wt\nu$ with $\wt \nu((0,t] \times E) = tJ(E \cap B(x,R_0)^c)$ for $E \in \BB(\R^d\setminus\{0\})$.
 In particular, by \eqref{e:C0-1} we have $\P(\wt Z_t = 0) \ge 1 - t J(B(0,R_0)^c) \ge 1-c_7 t$. Since $X^3$ and $\wt Z$ are independent, for any $x_0 \in \R^d$, $t>0$, $r \le R_0/2$, $x \in B(x_0,r)$, and a Borel set $F \subset B(x_0,r)$,
 $$ \P_x(X^{B(x_0,r)}_t \in F) \ge  \P_{x}(\{ X^{3,B(x_0,r)}_t \in F \}) \P(\wt Z_t = 0) \ge (1-c_7 t)\P_{x}(\{ X^{3,B(x_0,r)}_t \in F \}). $$

Now we choose $r_0 \in (0,R_0/2)$ small so that $1-c_7 \phi(\delta_3 r_0) \ge \frac{1}{2}$.
Then, using \eqref{e:1234}, we have for any $x_0 \in \R^d$, $r \le r_0$, and $0<t \le \phi(\delta_3 r)$,
\[ 
p^{B(x_0,r)}_{X}(t,x,y) \ge (1-c_7 t)p^{B(x_0,r)}_{X^3}(t,x,y) \ge \frac{c_6}{2}\phi^{-1}(t)^{-d}, \quad  x,y \in B(x_0,\delta_3 \phi^{-1}(t)).
\]
This finishes the proof.
\qed

Following the proof of \cite[Propositions 4.9 and 4.11]{CKL}, we obtain upper and lower bounds for the survival probability $\P_{x}(\tau_{B(x,r)}^X \le t)$ from \eqref{e:NDL}. 
\begin{lem}\label{l:survival}
	Let $X$ be a symmetric L\'evy process satisfying \eqref{e:TJ} and \eqref{e:NDL}, where $\psi$ satisfies \eqref{e:psi} with $\alpha \in (0,2]$ and $\phi$ is defined as \eqref{e:phi}. Then, there exist constants $c_3,c_4,c_5,c_6,r_0>0$ such that
	\begin{equation}\label{e:survival}
		\P_{x}(\tau_{B(x,r)}^X \le t) \le c_3 \left(\frac{t}{\psi(r/4)} + \exp \big( -c_4 \frac{r}{t^{1/2}} \big) \right) \quad  \mbox{for all} \quad t>0,\,\, r \le r_0, \,\, x \in \R^d,
	\end{equation}
	\begin{equation}\label{e:survival1}
	\P_{x}(\tau_{B(x,r)}^X \le t) \le \frac{c_5 t}{\phi(r)}\quad  \mbox{for all} \quad t>0,\,\,r \le r_0, \,\, x \in \R^d,
	\end{equation}
	and 
	\begin{equation}\label{e:survival2}
		\P_{x}(\tau_{B(x,r)}^X \le t) \ge \frac{c_6 t}{\psi(4r)}\quad  \mbox{for all} \quad 0<t \le \psi(4r)\,\, r \le r_0, \,\, x \in \R^d.
	\end{equation}
			Furthermore, if $A \neq 0$ is non-degenerate, we have
		\begin{equation}\label{e:survival3}
			\P_{x}(\tau_{B_{H}(x,r)}^{X^{H}} \le t) \ge c_8\exp\big(-c_7 \frac{r^2}{t}\big) \quad \mbox{for all } \,\, 0<t \le r^2,\,\, r \le r_0, \, x \in \R^d.
	\end{equation}		
The constants $c_i, i\in\{3,4,5,6,7,8\}$, and $r_0$ depend only on $\alpha$, $A$, $C_\psi$, $d$, the constants $c_1,c_2$ in \eqref{e:TJ}, and $c,\delta,r_0$ in \eqref{e:NDL}.
\end{lem}
\pf 
The conditions \eqref{e:TJ} and \eqref{e:NDL} imply  \cite[Assumption 4.4 and ${\rm TJ}_{1}(\psi,\le)$]{CKL} (also, see \cite[Proposition 4.3(i)]{CKL}). 
Now, the proof is similar to \cite[Propositions 4.9 and 4.11]{CKL}, and we only provide the differences. 
Note that by the definition of $\phi$ we have for $0<r \le R \le R_0$,
\[ 
\frac{\phi(R)}{\phi(r)} = \frac{R^2}{r^2} \cdot  \frac{\Vert A \Vert + 2\int_0^r \frac{sds}{\psi(s)}}{\Vert A \Vert + 2\int_0^R \frac{sds}{\psi(s)}} \le \frac{R^2}{r^2}, 
\]
so that the condition ${\rm US}_{R_0}(\phi,2,1)$ hold (see \cite[Definition 1.9]{CKL} for definition). 
Thus, the only difference is that we do not impose the condition ${\rm US}_{R_0}(\psi,\beta_2, C'_{U})$ in \cite[Proposition 4.9]{CKL}.
Since the condition ${\rm US}_{R_0}(\psi,\beta_2, C'_{U})$ has been used only for $\psi(r) \le C'_U 4^{\beta_2}\psi(r/4)$, we obtain that for any $t>0$, $r \le r_0$ and $x \in \R^d$,
\[ 	
\P_{x}(\tau_{B(x,r)}^X \le t) \le c_3 \left(\frac{t}{\psi(r/4)} + \exp \big( -c_4 \frac{r}{\vartheta_1(t,Cr)} \big) \right), 
\]
and
\[ \P_{x}(\tau_{B(x,r)}^X \le t) \le \frac{c_5 t}{\phi(r)}, \]
where the function $\vartheta_1(t,Cr) \in [\psi^{-1}(t),\phi^{-1}(t)]$ is defined in \cite[(4.2)]{CKL}. By \eqref{e:inverse}, we have $\vartheta_1(t,Cr) \le \phi^{-1}(t) \le ct^{1/2}$ and this implies \eqref{e:survival}.

For \eqref{e:survival2} and \eqref{e:survival3}, we follow the proof of \cite[Proposition 4.11(i) and 4.12(i)]{CKL}. Since the condition ${\rm TJ}_{R_0}(\psi(2\,\cdot), \ge)$ holds from \eqref{e:TJ}, it follows from the last inequality in \cite[Page 661]{CKL},
$$ \P_{x}(\tau_B \le t) \ge 1 - e^{-c t / \psi(4r)} \ge \frac{c_6 t}{\psi(4r)}, \quad 0<t \le \psi(4r). $$
This proves \eqref{e:survival2}. Note that since $\phi(r) \asymp r^2$ when $A \neq 0$ is non-degenerate, we can apply \cite[(4.38)]{CKL} by \eqref{e:NDL}. Thus,
\begin{equation}\label{e:survival31}
 \P_x(\tau_B \le t) \ge \exp \big( -c_7\frac{r}{\vartheta_2(t,Cr)} \big), 
\end{equation}
where the function $\vartheta_2(t,Cr) \in [\psi^{-1}(t), \phi^{-1}(t)]$ is defined in
 \cite[(4.3)]{CKL}. To obtain proper lower bound of $\vartheta_2(t,r)$, observe the definition in \cite[(4.3)]{CKL} that if a number $\rho_0 \le \phi^{-1}(t)$ satisfies
 $ \frac{t\rho_0}{\phi(\rho_0)} \ge  r$,
 then $\vartheta_2(t,r) \ge \rho_0$, whether $\rho_0 \in [\psi^{-1}(t), \phi^{-1}(t)]$ or $\rho_0 < \psi^{-1}(t)$. By \eqref{e:inverse} and $t \le r^2 \le r_0^2$ we have
 $$ \Vert A \Vert^{1/2} \frac{t}{r} \le \Vert A \Vert^{1/2}t^{1/2} \le \phi^{-1}(t). $$
 Taking $\rho_0 := \big( \frac{\Vert A \Vert}{C} \land \Vert A\Vert^{1/2} \big) \frac{t}{r} \le \phi^{-1}(t)$, by \eqref{e:square} we have
 $         \frac{t\rho_0}{\phi(\rho_0)} \ge \Vert A \Vert \frac{t}{\rho_0} \ge Cr, $
 which implies $\vartheta_2(t,Cr) \ge \rho_0$. This combined with \eqref{e:survival31} yields \eqref{e:survival3}.

 \qed

	\begin{cor}\label{c:5.3}
		If \eqref{e:survival1} holds, then there exists a constant $c>0$ such that
		$$ \E_x[\tau_{B(x,r)}^X] \ge \frac{\phi(r)}{4c_5}. $$
	\end{cor}
	\pf Taking $t=\frac{\phi(r)}{2c_5}$ at \eqref{e:survival1} we have
	\[ 
	\E_x[\tau_{B(x,r)}^X] \ge \frac{\phi(r)}{2c_5} \P_x(\tau_{B(x,r)}^X > \frac{\phi(r)}{2c_5}) =  \frac{\phi(r)}{2c_5} \big( 1- \P_x(\tau_{B(x,r)}^X \le \frac{\phi(r)}{2c_5}) \big) \ge \frac{\phi(r)}{4c_5}. 
	\]
	\qed

\subsection{Projections of $X$ onto Subspaces}\label{subsection:projection}
	In this subsection, we study the projection of the L\'evy process $X$. 
	Let $\{0\} \neq H \subset \R^d$ be a subspace and $d_1 := {\rm dim} \, H$. Let $\pi_H : \R^d \to H$ be the projection from $\R^d$ to $H$. We use $B_{H}(x,r):=B(x,r)\cap H$ to denote an open ball at $x$ with radius $r>0$ in $H$.
Let $\pi_H : \R^d \to H$ be the projection from $\R^d$ to $H$. The process $X^H = \{X^H_t\}_{t \ge 0}$ defined by $X^H_t := \pi_H(X_t)$ is a $d_1$-dimensional L\'evy process in $H$, where $d_1= {\rm dim} \, H$ is the vector space dimension of the subspace. Let $J_H(dx)$ be the L\'evy measure of $X^H$.
	\begin{lem}\label{l:projection}
		 Let $X$ be a symmetric L\'evy process satisfying \eqref{e:J} and \eqref{e:psi}. Also, assume that $A=0$ or $A$ is non-degenerate. Then, 
		
		\noindent {\rm (i)} There exist constants $c_1,c_2>0$ and $r_1>0$ such that for any subspace $\{0\} \neq H \subset \R^d$ and $r \le r_1$,
		$$ 
		\frac{c_1}{\psi(2\sqrt{d}r)} \le J_H(B_{H}(0,r)^c) \le \frac{c_2}{\psi(r)}, 
		$$
		where $B_{H}(0,r) = \{ y \in H : |y| < r \}$ is an open ball in $H$.
		
		\noindent {\rm (ii)} 
		There exist $c_3>0$, $r_0>0$ and $\delta_1 \in (0,1)$ such that for any for any subspace $\{0\} \neq H \subset \R^d$ and $x_0 \in H$,
		\begin{equation}\label{e:NDL1} 
		p_{X^H}^{B_{H}(x_0,r)}(t,x,y) \ge c_3\phi^{-1}(t)^{-d_1}, \qquad r \le r_0, \,\, 0<t \le \phi(\delta_1 r), \,\, x,y \in B_{H}(x_0,\delta_1 \phi^{-1}(t)),
		\end{equation}
			where $r_0$ and $\delta$ are the constants in \eqref{e:NDL}, and the function $\phi$ is defined in \eqref{e:phi}.
		The constants $c_1,c_2,c_3,r_0$ and $\delta_1$ depend only on $\alpha$, $A$, $C_1$, $C_2$, $C_3$, $C_\psi$, $d$ and $R_0$, and are independent of $H$.
	\end{lem}
	\pf 
First, we establish $\rm{(i)}$. Without loss of generality, we may assume that $H = \{(x_1, \dots, x_{d_1},0,\dots, 0) : x_i \in \R, 1 \le i \le d_1 \}$. Note that $d_1 \ge 1$ since $\{0\} \neq H$. Observe that
	\begin{align*}
		J_H(B_H(0,r)^c) = J(\pi_H^{-1}(B_H(0,r)^c)) \le J(B(0,r)^c).
	\end{align*}
Thus, by Lemma \ref{l:TJ} we have $J_H(B_H(0,r)^c) \le \frac{c_1}{\psi(r)}$ for any $r \le R_0/2$. 

For the lower bound, consider
	$$ C:= \{ (x_1, \dots, x_d) : r \le x_i \le 2r, \,\, 1 \le i \le d   \} \subset \pi_H^{-1}(B_H(0,r)^c). $$ 
	Note that since
	$$ \pi_H(C) = \{ (x_1,\dots,x_{d_1},0,0,\dots,0) : r \le x_i \le 2r, \,\, 1 \le i \le d_1 \} \subset B_H(0,r)^c, $$
We have 
	\begin{align*}
	J_H(B_H(0,r)^c) &\ge J(C) = \int_r^{2r} \dots  \int_r^{2r} \int_{r}^{2r} \frac{C_1}{|x|^d \psi(|x|)} \ud x_1 \dots \ud x_d \\
		&\ge  \frac{C_1 r^d }{(2\sqrt d r)^d \psi(2 \sqrt d r)} \ge \frac{c_2}{\psi(2 \sqrt d r)}.
	\end{align*}

Now we prove $\rm{(ii)}$. Let $r_0$ and $\delta$ be the constants in \eqref{e:NDL}. For any $r \le r_0$, $0<t \le \phi(\delta r)$, $x^0=(x^0_1, \dots, x^0_d) \in H$, $x \in B_H(x^0,\delta \phi^{-1}(t)/\sqrt d )$ and any Borel set $A \subset B_H(x^0,\delta \phi^{-1}(t)/ \sqrt d)$, we have
\begin{align*}
&\int_A p^{B_H(x^0,r)}_{X^H}(t,x,y)\ud y = \P_x(X^{H,B_H(x^0,r)}_t \in A) \\ =& \P_x(X_t^{\pi_H^{-1}(B_H(x^0,r))} \in \pi_H^{-1}(A)) \ge \P_x(X_t^{B(x^0,r)} \in \pi_H^{-1}(A)), 
\end{align*}
where $X^{H,B_H(x^0,r)}$ is the killed process of $X^H$ upon $B_H(x^0,r)$.
Let $F := \prod_{i=d_1 + 1}^d B_{\R}(x^0_i, \delta \phi^{-1}(t) / \sqrt d)$ be a box centered in $(x_{d_1+1}^{0},\cdots,x_{d}^{0}) \in  \R^{d-d_1}$. Then, for any $A \subset H$ we have $A \times F \subset \pi_H^{-1}(A) \cap B(x^0,\delta \phi^{-1}(t))$. Thus by \eqref{e:NDL} we have 
\begin{align*}
 \P_x(X_t^{B(x^0,r)} \in \pi_H^{-1}(A)) &\ge \P_x(X_t^{B(x^0,r)} \in A \times F) = \int_{A \times F} p_X^{B(x^0,r)}(t,x,y)\ud y \\ &\ge c_3\phi^{-1}(t)^{-d} |A \times F| = c_4 \phi^{-1}(t)^{-d_1} |A|.   \end{align*}
Define $\delta_1=\delta/\sqrt{d}$. Then, since $\delta_1 < \delta$, \eqref{e:NDL1} follows from the above estimate, the Lebesgue differentiation theorem, and the joint continuity of the Dirichlet heat kernel.

Note that every constant in this proof is independent of nontrivial subspace $H$.
\qed

\noindent From Lemmas \ref{l:survival} and \ref{l:projection}, we conclude that
\begin{cor}\label{c:sur}
		 Let $X$ be a symmetric L\'evy process satisfying \eqref{e:J} and \eqref{e:psi}. Also, assume that $A=0$ or $A$ is non-degenerate. Then, there exist constants $c_i$, $i\in\{1,2,3,4,5\}$, such that for any subspace $\{0\} \neq H \subset \R^d$,
		\begin{equation}\label{e:sur}
		\P_{x}(\tau_{B_{H}(x,r)}^{X^{H}} \le t) \le c_1 \left(\frac{t}{\psi(r/4)} + \exp \big( -c_2 \frac{r}{t^{1/2}} \big) \right) \quad \mbox{for all } \,\, t>0,\,\, r \le r_0, \, x \in \R^d,
		\end{equation}
		\begin{equation}\label{e:sur1}
		\P_{x}(\tau_{B_{H}(x,r)}^{X^{H}} \le t) \le \frac{c_3 t}{\phi(r)} \quad \mbox{for all } \,\, t>0, \,\, r \le r_0, \, x \in \R^d,
		\end{equation}
		and 
		\begin{equation}\label{e:sur2}
			\P_{x}(\tau_{B_{H}(x,r)}^{X^{H}} \le t) \ge \frac{c_4 t}{\psi(4 \sqrt{d} r)} \quad \mbox{for all } \,\, 0<t \le \psi(4 \sqrt{d} r),\,\, r \le r_0, \, x \in \R^d.
		\end{equation}
		Furthermore, if $A \neq 0$ is non-degenerate, we have
		\begin{equation}\label{e:sur3}
		\P_{x}(\tau_{B_{H}(x,r)}^{X^{H}} \le t) \ge c_6\exp\big(-c_5 \frac{r^2}{t}\big) \quad \mbox{for all } \,\, 0<t \le r^2,\,\, r \le r_0, \, x \in \R^d.
	\end{equation}		
The constants $c_i, i\in\{1,2,3,4,5,6\}$, and $r_0$ depend only on $\alpha$, $A$, $C_1$, $C_2$, $C_3$, $C_\psi$ and $d$.
 \end{cor}

\vspace{0.2in}
\noindent
\textbf{Acknowledgments:} 
Part of this project was conducted while the second-named author was visiting the Korea Institute for Advanced Study (KIAS) in 2024. He expresses his gratitude to KIAS for its support, warm hospitality, and excellent working environment.

\vskip 0.3truein

{\bf Jaehun Lee}

Stochastic Analysis and Applicaton Research Center, KAIST,
Daejeon 34141, South Korea

E-mail: \texttt{hun618@kaist.re.kr}

\bigskip

{\bf Hyunchul Park}

Department of Mathematics, State University of New York at New Paltz, NY 12561,
USA

E-mail: \texttt{parkh@newpaltz.edu}

\end{document}